\documentclass{article}

\usepackage{amsmath}
\usepackage{amssymb}
\usepackage{latexsym}
\usepackage{amsxtra}
\usepackage{amscd}
\usepackage{theorem}
\usepackage{xypic}
\usepackage{verbatim}

\theoremstyle{plain}

{\theorembodyfont{\slshape}

        \newtheorem{thm}{Theorem}[section]
        \newtheorem{cor}[thm]{Corollary}
        \newtheorem{lem}[thm]{Lemma}
        \newtheorem{prop}[thm]{Proposition}

}

{\theorembodyfont{\rmfamily}

        \newtheorem{defn}[thm]{Definition}
        \newtheorem{prob}[thm]{Problem}
        \newtheorem{rem}[thm]{Remark}
        
        \newtheorem{exa}[thm]{Example}
        \newtheorem{ass}[thm]{Assumption}
        
        \newtheorem{claim}[thm]{Claim}

        \newtheorem{conj}[thm]{Conjecture}
}

\renewcommand{\em}{\sl}

\newcommand{\proof}{{\bf Proof:\ }}
\newcommand{\Endproof}{\hspace*{\fill} $\Box$ \vspace{1ex} \noindent }

\makeatletter
\renewcommand{\subsection}{\@startsection{subsection}{2}%
        {\z@}{-3.25ex plus -1ex minus-.2ex}{-1em}{\bf}}
\makeatother

\newcommand{\NN}{\mathbb{N}}
\newcommand{\ZZ}{\mathbb{Z}}
\newcommand{\QQ}{\mathbb{Q}}
\newcommand{\FF}{\mathbb{F}}
\newcommand{\RR}{\mathbb{R}}

\newcommand{\PP}{\mathbb{P}}
\newcommand{\KK}{\mathbb{K}}
\newcommand{\GG}{\mathbb{G}}

\renewcommand{\AA}{\mathbb{A}}
\newcommand{\BB}{\mathbb{B}}

\newcommand{\C}{\mathcal{C}}
\newcommand{\G}{\mathcal{G}}

\newcommand{\OO}{\mathcal{O}}

\newcommand{\K}{\mathfrak{K}}

\newcommand{\Gal}{{\rm Gal}}

\newcommand{\m}{\mathfrak{m}}

\newcommand{\Spec}{{\rm Spec\,}}

\newcommand{\Spm}{{\rm Spm}}
\newcommand{\Frac}{{\rm Frac}}

\newcommand{\ord}{{\rm ord}}

\newcommand{\inj}{\hookrightarrow}

\newcommand{\abs}[1]{\lvert#1\rvert}

\newcommand{\sw}{{\rm sw}}
\newcommand{\dsw}{{\rm dsw}}

\newcommand{\ol}[1]{\overline{#1}}
\newcommand{\blockmatrix}[4]{
    \left( \begin{array}{c|c} #1 & #2 \\ \hline #3 &
    #4 \end{array} \right)}

\newcommand{\Xb}{\bar{X}}
\newcommand{\Yb}{\bar{Y}}
\newcommand{\Zb}{\bar{Z}}
\newcommand{\Wb}{\bar{W}}
\newcommand{\Db}{\bar{D}}

\newcommand{\Ab}{\bar{A}}

\newcommand{\Kb}{\bar{K}}

\newcommand{\xb}{\bar{x}}
\newcommand{\yb}{\bar{y}}

\newcommand{\fb}{\bar{f}}

\newcommand{\cb}{\bar{c}}
\newcommand{\chib}{\bar{\chi}}

\newcommand{\KKh}{\hat{\KK}}

\makeatletter
\newcommand{\subjclass}[2][2010]{%
  \let\@oldtitle\@title%
  \gdef\@title{\@oldtitle\footnotetext{#1 \emph{Mathematics subject classification.} #2}}%
}
\newcommand{\keywords}[1]{%
  \let\@@oldtitle\@title%
  \gdef\@title{\@@oldtitle\footnotetext{\emph{Key words and phrases.} #1.}}%
}
\makeatother
\begin{document}

\title{Cyclic extensions and the local lifting problem}
\author{Andrew Obus and Stefan Wewers}

\maketitle

\begin{abstract}
The local Oort conjecture states that, 
if $G$ is cyclic and $k$ is an algebraically closed field of characteristic $p$, 
then all $G$-extensions of $k[[t]]$ should lift to characteristic zero.  We prove a critical case of this conjecture.  In particular, 
we show that the conjecture is always true when $v_p(|G|) \leq 3$, and is true for arbitrarily highly $p$-divisible 
cyclic groups $G$ when a certain condition on the higher ramification filtration is satisfied.  
\end{abstract}

\section{Introduction}\label{Sintro}
\subsection{The local lifting problem}\label{Sllp}
Let $\Yb$ be a smooth, projective, connected curve over an algebraically closed field $k$ of characteristic $p>0$.  Results in 
deformation theory going back to Grothendieck show that $\Yb$ can always be lifted to characteristic zero.
Specifically, one can always find a discrete valuation ring (DVR) $R$ in characteristic zero, with residue 
field $k$, such that there exists a smooth relative $R$-curve $Y$ with $Y \times_R k \cong \Yb$.  

In \cite{Oort87}, Oort asked the natural question: can one lift a \emph{Galois cover} of curves to characteristic zero?  That is, if $G$ is a finite group, and 
$\fb:\Yb \to \Xb$ is a $G$-Galois cover of smooth, projective, connected curves, is there a $G$-Galois cover $f: Y \to X$ of 
smooth relative curves over a DVR $R$ in characteristic zero whose special fiber is $\fb: \Yb \to \Xb$?  Clearly, the answer is not 
always ``yes."  For instance, the group $G = \ZZ/p \times \ZZ/p$ acts faithfully on $\Yb = \PP^1_k$ via an embedding into
$\GG_a(k)$.
If $\Xb = \Yb/G$, then the generic fiber of any lift $f: Y \to X$ of $\fb : \Yb \to \Xb$ must be a $G$-Galois cover $\PP^1_K \to \PP^1_K$,
where $K = \Frac(R)$.  But 
if $p \geq 3$, then $G$ cannot act faithfully on $\PP^1$ in characteristic zero, so such a lift does not exist.  
However, Oort conjectured  (\cite{Oort95}) that $G$-covers
should always lift when $G$ is cyclic.  Our main result (Theorem \ref{Tmain}) proves 
the Oort conjecture whenever $v_p(|G|) \leq 3$, and in many cases for arbitrarily large cyclic groups $G$.  
Specific statements are in \S\ref{Snew}.  The case we prove is critical, as Pop (\cite{Pop}) is able to reduce
the conjecture to the case we have proven, thus proving the entire conjecture.  See Remark \ref{Roort}.  

Although this was not yet known at the time of \cite{Oort87}, it turns out that the nature of lifting Galois covers of curves is completely local.  Let $B$
be the branch locus of the $G$-Galois cover $\fb: \Yb \to \Xb$.  If, for each $x \in B$, one can lift the cover 
$\fb$ when restricted to a formal neighborhood $\hat{x}$ of $x$, then one can lift $f$ as well.  
This is known as the \emph{local-global principle}.  Proofs have been given by
Bertin and M\'{e}zard (\cite{BM}), Green and Matignon (\cite{GM98}), and Garuti (\cite{Garuti96}).  

The restriction $\fb|_{\hat{x}}$ is a disjoint union of covers of the form $\Spec k[[z]] \to \Spec k[[t]]$.  Thus, the study of lifting Galois
covers can be reduced to the following \emph{local lifting problem}:

\begin{prob}[The local lifting problem]\label{llp}
  Let $k$ be an algebraically closed field of characteristic $p$ and $G$ a finite group.  
  Let $k[[z]]/k[[t]]$ be a $G$-Galois extension (i.e., an integral extension of integrally closed domains that is $G$-Galois on the level of fraction fields).
  Does there exist a DVR $R$ of characteristic zero with residue field $k$ and a
  $G$-Galois extension $R[[Z]]/R[[T]]$ that reduces to $k[[z]]/k[[t]]$?  That is, does the $G$-action on $R[[Z]]$ reduce to that on $k[[z]]$,
  if we assume that $Z$ (resp.\ $T$) reduces to $z$ (resp.\ $t$)?
\end{prob}

\begin{rem}\label{ringfieldrem}
If $L/k[[t]]$ is any $G$-Galois extension,
then the structure theorem for complete DVRs shows that $L$ is always abstractly isomorphic to $k[[z]]$.  
So in the local lifting problem, we may as well talk about general $G$-Galois extensions of $k[[t]]$.
\end{rem} 

In light of the local-global principle, one translates the conjecture of Oort above into the local context. 

\begin{conj}[Local Oort conjecture] \label{OortConj} 
  The local lifting problem can always be solved when $G$ is cyclic.
\end{conj} 

Note that a consequence of the local Oort conjecture and the local-global principle is that
any Galois cover of $k$-curves with cyclic \emph{inertia} groups (not just cyclic Galois group) lifts to characteristic zero.  

The first author's paper \cite{Ob11} is a detailed exposition of many aspects of the local lifting problem.

\subsection{Previously known results}\label{Sknown}
	It is easy to prove that if Conjecture \ref{OortConj} is true for a given
$\ZZ/p^n$-extension $L_n/k[[t]]$, then it is also true for any
$\ZZ/rp^n$-extension ($p \nmid r$) whose unique
$\ZZ/p^n$-subextension is $L_n/k[[t]]$ (see, e.g., the proof of \cite{Ob11},
Proposition 6.3).  Thus, we may reduce to the case where $G \cong \ZZ/p^n$. 

If $G \cong \ZZ/p$ or $\ZZ/p^2$, then Conjecture \ref{OortConj} has already been shown to be true.  
The original proof of the case $G \cong \ZZ/p$ was given by Oort, Sekiguchi, and Suwa in \cite{OSS}.  Their proof was given in the global
context, for any $\ZZ/p$-cover of smooth, projective curves.  The proof used essentially global techniques, such as generalized Jacobians and 
abelian schemes.  

Later, Green and Matignon reproved Conjecture \ref{OortConj} for $G = \ZZ/p$, and also proved it for $G=\ZZ/p^2$ 
(\cite{GM98}, II, Theorems 4.1 and 5.5).  Their lifts were given by explicit Kummer extensions.  The form of
these extensions was inspired by \emph{Sekiguchi-Suwa Theory} (or \emph{Kummer-Artin-Schreier-Witt Theory}), although 
Green and Matignon developed their proofs independently.  This theory
gives (in principle) explicit equations for group schemes classifying unramified $\ZZ/p^n$-extensions of flat local $R$-algebras, 
where $R$ is a complete discrete valuation ring in mixed characteristic $(0, p)$.  When $n \leq 2$, it is manageable to write down
these equations explicitly, and Green and Matignon were able to exploit this to write down their lifts.  
See \cite{SS94} for an overview of the general theory, \cite{SS99} for an expanded version with proofs included, 
\cite{SS01} for a detailed account of the case $n \leq 2$, or \cite{Tossici10} for a briefer overview of this case. 

Unfortunately, when $n \geq 3$, the equations involved in Sekiguchi-Suwa theory become extremely complicated,
and extraordinarily difficult to work with.  No one has been able to use the method of \cite{GM98} to prove Conjecture \ref{OortConj}
for \emph{any} $\ZZ/p^n$-extension with $n \geq 3$.  Indeed, prior to this paper, Conjecture \ref{OortConj} 
for such extensions was only known to be true 
for sporadic examples arising from Lubin-Tate formal group laws (\cite{GM99}, \cite{Gre04}).

However, there has long been evidence for the truth of Conjecture \ref{OortConj}, in the sense that all of the main known obstructions to lifting (such as the 
Bertin/Katz-Gabber-Bertin obstructions of \cite{Bertin}, \cite{CGH1}, and the Hurwitz tree obstruction of \cite{Brewis}, \cite{Brewisthesis}) vanish for cyclic
extensions.  

\subsection{Main result}\label{Snew}
To state our main result, Theorem \ref{Tmain}, we recall that a $\ZZ/p^n$-extension $L_n/k[[t]]$ 
gives rise to a higher ramification filtration $G^s_{s \geq 0}$ on the group $G$ for the 
upper numbering (\cite{SerreCL}, IV).  The breaks in this filtration (i.e., the values $i$ for which $G^i \supsetneq G^j$ for all $j > i$) will
be denoted by $(m_1, m_2, \ldots, m_n)$.  One knows that $m_i \in \NN$ and
\begin{equation} \label{miineq}
  m_i\geq p\,m_{i-1},
\end{equation}
for $i=2,\ldots,n$ (see, e.g., \cite{Garuti99}). 

\begin{thm}\label{Tmain}
Let $L_n/k[[t]]$ be a $\ZZ/p^n$-extension with upper ramification breaks 
$(m_1, m_2, \ldots, m_n)$.  Suppose, for $3 \leq i \leq n-1$, that there is never an integer $a_i$
such that $\frac{m_i}{p} - m_{i-1} < a_i \leq \left(\frac{m_i}{m_i - m_{i-1}}\right)\left(\frac{m_i}{p} - m_{i-1}\right)$.  
Then Conjecture \ref{OortConj} holds for $L_n/k[[t]]$.
\end{thm}

\begin{rem}\label{Roort}
\begin{enumerate}
\item It is not hard to see that the condition in Theorem \ref{Tmain} 
is satisfied whenever $L_n/k[[t]]$ has no ``essential ramification,"  i.e., that $m_i < pm_{i-1} + p$ for $2 \leq i \leq n$.
Pop's proof of the Oort conjecture in  in \cite{Pop} reduces the (local) Oort conjecture to the case where
there is no essential ramification.
\item The condition in Theorem \ref{Tmain} is vacuous for $n = 3$, so Conjecture \ref{OortConj} for $G \cong \ZZ/p^3$ is an immediate consequence.
\item The condition that there is no integer $a_i$ such that 
$\frac{m_i}{p} - m_{i-1} < a_i \leq \left(\frac{m_i}{m_i - m_{i-1}}\right)\left(\frac{m_i}{p} - m_{i-1}\right)$ is equivalent to saying that, if 
$m_i = pm_{i-1} + pr_i - \eta_i$, with $r_i$ and $\eta_i$ integers such that $0 \leq \eta_i < p$, then $0 \leq r_i \leq \eta_i$.  See \cite{Pop}.
\item One can conjecture further that, if $G \cong \ZZ/p^n$, one should be able to take $R = W(k)[\zeta_{p^n}]$ in Conjecture \ref{OortConj},
where $\zeta_{p^n}$ is a $p^n$th root of unity.
This is known when $n \leq 2$.  Unfortunately, our proof gives no effective bounds on $R$.  
The proof of \cite{Pop} gives, in theory, some effective bounds on $R$, but they are much weaker than what is expected.
\end{enumerate}
\end{rem}

\begin{rem}\label{meta}
It would be interesting to investigate the local lifting problem when $G \cong \ZZ/p^n \rtimes \ZZ/m$, with $p \nmid m$.  One can show
that a necessary condition to lift a $G$-extension $L/k[[t]]$ to characteristic zero is that the action of $\ZZ/m$ on $\ZZ/p^n$ is either faithful 
or trivial, and that, if it is faithful, then the upper ramification breaks $(m_1, \ldots, m_n)$ of the $\ZZ/p^n$-subextension are all congruent to
$-1 \pmod{m}$.  In light of \cite{Pop}, one can ask if this is the only restriction.  For instance, should all $D_{p^n}$-extensions lift for odd $p$?
\end{rem}

\subsection{Outline of the paper}
We start with a short section (\S\ref{Switt}) overviewing the 
basics of the Artin-Schreier-Witt theory, giving an explicit characterization of the $\ZZ/p^n$-extensions of $k[[t]]$.
In \S\ref{Sind}, we set up our induction on $n$, and show how it proves Theorem \ref{Tmain}. 
We prove the base cases $n = 1$ and $n=2$ in \S\ref{plan2}.  The paper begins
in earnest with \S\ref{Sgeom}.  In \S\ref{Sgeomsetup} and \S\ref{Scharacters} we introduce the language of \emph{characters}, which 
will often be more convenient than the language of extensions for expressing our results.  In \S\ref{swan}, we introduce Kato's  
\emph{Swan conductor} in the situation relevant to us.  The Swan conductor serves several purposes in this paper, most notably
giving us a way to measure how
bad the reduction of a cover is.  In \S\ref{orderp} and \S\ref{Skink}, we examine the particular case of $\ZZ/p$-extensions in great
detail.  This is important, as $\ZZ/p$-extensions are the building blocks of our inductive process.

In \S\ref{Sproof}, we give the main proofs.  Unlike in \cite{GM98}, we do not try to write down a lift of a given $\ZZ/p^n$-extension explicitly.
In particular, we do not use the Sekiguchi-Suwa theory at all, except in the relatively trivial case of $\ZZ/p$-extensions
(i.e., Kummer-Artin-Schreier theory).
Instead, we write down what the form of the equations should be, in order that we might lift
\emph{some} $\ZZ/p^n$-extension.  
Then, we show that if equations in this form do not reduce to a Galois extension,
then they can be deformed to yield something that comes closer to reducing to a Galois extension, and 
this deformation process eventually terminates.  This proves Conjecture \ref{OortConj} for some $\ZZ/p^n$-extension.  
We then show that, given a solution for some particular extension, we can find solutions to many more.
A more detailed outline of \S\ref{Sproof} is given in \S\ref{plan}.  

The proofs of several key technical results are postponed to \S\ref{Stechnical}.  This is partially because
proving these results would disrupt the continuity of \S\ref{Sproof}, and partially because the proofs share much notation, and thus
are more easily read together.

\subsection{Conventions}
The letter $K$ will always be a field of characteristic zero that is complete with respect to a
discrete valuation $v:K^\times\to\QQ$. We assume that the residue field $k$ of
$K$ is algebraically closed of characteristic $p$. We also assume that the
valuation $v$ is normalized such that $v(p)=1$.  The ring of
integers of $K$ will be denoted $R$.  We fix an algebraic closure $\Kb$ of $K$, and whenever necessary, we will replace $K$ by a suitable finite
extension within $\Kb$, without changing the above notation.  The maximal ideal of $R$ will be denoted $\m$.
Furthermore, for each $r \in \NN$, we fix 
once and for all a compatible system of $r$th roots $p^{1/r}$ in $\Kb$ such that
if $ab = r$, then $(p^{1/r})^a = p^{1/b}$.  Thus $p^q \in \Kb$ is well defined for any $q \in \QQ$.
The greatest integer function of $x$ is written $[x]$.
%
%
%
\section*{Acknowledgements}
The authors would like the thank Irene Bouw, Louis Brewis, David Harbater, and Martin Rubey for useful conversations.
In particular, Louis Brewis showed us that the Hurwitz tree obstruction to the local Oort conjecture vanishes, which helped us come up with the
idea of where to place our branch points.  Irene Bouw showed us how to vastly simplify our proof of Theorem \ref{thm2}.
We also thank Michel Matignon and Florian Pop for helpful comments on earlier versions of this article.  

The first author was partially supported by an NSF Mathematical Sciences Postdoctoral Research Fellowship.  This work
was completed during a visit by the first author to the Max-Planck-Institut f\"{u}r Mathematik, and he thanks the Institute for its support
and pleasant atmosphere.

\section{Artin-Schreier-Witt theory}\label{Switt}

If $L_n/k[[t]]$ is a $\ZZ/p^n$-extension, then so is the extension $M_n/k((t))$, where $M_n = \Frac(L_n)$.
The classical Artin-Schreier-Witt theory states that $M_n/k((t))$ is given by
an Artin-Schreier-Witt equation 
\[
       \wp(x_1,\ldots,x_n)=(f_1,\ldots,f_n),
\]
where $(f_1, \ldots, f_n)$ lies in the ring  $W_n(k((t)))$ of truncated Witt vectors, 
$F$ is the Frobenius morphism on $W_n(k((t)))$, and $\wp(x):=F(x)-x$ is the Artin-Schreier-Witt 
isogeny.  Adding a truncated Witt vector of the form $\wp(y)$ to $(f_1,\ldots,f_n)$ does not change the extension, 
and we obtain a group isomorphism $H^1(k((t)), \ZZ/p^n) \cong W_n(k((t)))/\wp(W_n(k((t))))$. 
Since we can add $\wp(y)$ to a Witt vector without changing the extension, we may assume that the 
$f_i$ are polynomials in $t^{-1},$ all of whose terms have prime-to-$p$ degree.  In this case, if
\begin{equation}\label{rambreaks}
     m_i :=\max\{\, p^{i-j}\deg_{t^{-1}}(f_j) \mid j=1,\ldots,i\,\},
\end{equation}
then the $m_i$ are exactly the breaks in the higher ramification filtration of $M_n/k((t))$ (\cite{Garuti99}, Theorem 1.1).  
From this, one easily sees that $p \nmid m_1$, that $m_i \geq pm_{i-1}$ for $2 \leq i \leq n$, and that if $p | m_i$, then $m_i = pm_{i-1}$.

For more details, see \cite{Witt37} or the exercises on page 330 of \cite{Lang}.

\section{The induction process}\label{Sind}
Let $L_n/k[[t]]$ be a $\ZZ/p^n$-extension.  
A theorem of Harbater-Katz-Gabber (\cite{Ha80}, \cite{Ka86}) shows that (after possibly changing the uniformizer $t$ of $k[[t]]$)
there exists a unique cover $\Yb_n \to \Xb := \PP^1_k$ that is \'{e}tale
outside $t=0$, totally ramified above $t=0$, and such that the formal completion of $\Yb_n \to \Xb$ at $t=0$ yields the extension $L_n/k[[t]]$.
The local-global principle thus shows that solvability of the local lifting problem from $L_n/k[[t]]$ is equivalent
to the following claim, which will be more convenient to work with:

\begin{claim} \label{claim1}
  Given a $G$-Galois extension $L_n/k[[t]]$, with $G \cong \ZZ/p^n$, then after possibly changing the uniformizer $t$ of $k[[t]]$,
  there exists a $G$-Galois cover $Y_n\to X:=\PP^1_K$ 
  (where $K$ is the fraction field of some DVR $R$ as above) with the following properties:
  \begin{enumerate}
  \item The cover $Y_n\to X$ has good reduction with respect to the standard
    model $\PP_R^1$ of $X$ and reduces to a $G$-Galois cover
    $\bar{Y}_n\to\bar{X}=\PP^1_k$ which is totally ramified above the point
    $t=0$ and \'etale everywhere else.
  \item
     The completion of $\bar{Y}_n\to\bar{X}$ at $t=0$ yields $L_n/k[[t]]$.
  \end{enumerate}
\end{claim}  

\begin{rem}
  Let $T$ be a coordinate of $\PP^1_R$ reducing to $t$.  Then Condition (i) and (ii) in Claim \ref{claim1} can be reformulated as
  follows:
\begin{enumerate}
\item
  The cover $Y_n\to X$ is \'etale outside the open disk
    \[
            D:= \{\, T \mid \abs{T} < 1 \,\}.
    \]
\item
  The inverse image of $D$ in $Y_n$ is an open disk.
\item
  If $A = R[[T]]\{T^{-1}\}$ is the ring 
  $$\left\{ \sum_{j \in \ZZ} a_jT^j \left. \right| a_j \in R,\ a_j \to 0 \text{ as } j \to - \infty \right\},$$ then 
  the cover $Y_n\to X$ is unramified when base changed to $\Spec A$ 
  (which corresponds to the ``boundary" of the disk $D$).  The
  extension of residue fields is isomorphic to the extension of fraction fields
  coming from $L_n/k[[t_0]]$.
\end{enumerate}
\end{rem}

If $R$ is a characteristic zero DVR with residue field $k$ and fraction field $K$, set  
$D(r) = \{ T \in \Kb \mid \; \abs{T} < \abs{p}^r\}$, using the non-archimedean absolute value on $K$ 
induced from the valuation.  

We prove Theorem \ref{Tmain} (in the context of Claim \ref{claim1}) by induction using the following base case (Lemma \ref{Lbase})
and induction step (Theorem \ref{Tsetup}).

\begin{lem}\label{Lbase}
If $n = 1$ (resp.\ $n=2$), let $L_n/k[[t]]$ be a $\ZZ/p^n$-extension with upper ramification break $m_1$ (resp.\ breaks $(m_1, m_2)$). 
Then there exists a $\ZZ/p^n$-cover $Y_n \to X$ satisfying Claim \ref{claim1} for $L_n/k[[t]]$, which is \'{e}tale outside
the open disk $D(r_n)$, where $r_n = \frac{1}{m_n(p-1)}$.
\end{lem}

\begin{thm}\label{Tsetup}
\begin{enumerate}
\item Suppose $n > 1$, and let $L_n/k[[t]]$ be a $\ZZ/p^n$-extension with upper ramification breaks $(m_1, \ldots, m_n)$ and
$\ZZ/p^{n-1}$-subextension $L_{n-1}/k[[t]]$.  
Suppose there exists a $\ZZ/p^{n-1}$-cover $Y_{n-1} \to X$  satisfying Claim \ref{claim1} for $L_{n-1}/k[[t]]$, which is furthermore
\'{e}tale outside the open disk $D(r_{n-1})$, where $r_{n-1} = \frac{1}{m_{n-1}(p-1)}$.
Then there is a $\ZZ/p^n$-cover $Y_n \to X$ satisfying Claim \ref{claim1} for $L_n/k[[t]]$.
\item If $L_n/k[[t]]$ is as in part (i) and there is no integer $a$ satisfying
$$\frac{m_n}{p} - m_{n-1} < a \leq \left(\frac{m_n}{m_n - m_{n-1}}\right) \left(\frac{m_n}{p} - m_{n-1}\right),$$ then the 
$\ZZ/p^n$-cover
$Y_n \to X$ in part (i) can be chosen to be \'{e}tale outside $D(r_n)$, where $r_n = \frac{1}{m_n(p-1)}$.
\end{enumerate}
\end{thm}

Theorem \ref{Tmain} now follows easily.\\

\proof (of Theorem \ref{Tmain}) 
Let $L_n/k[[t]]$ be in the form of Theorem \ref{Tmain}. 
We note that, for any $n' \leq n$, the unique $\ZZ/p^{n'}$-subextension $L_{n'}/k[[t]]$ of $L_n/k[[t]]$
has upper ramification breaks $(m_1, \ldots, m_{n'})$ (\cite{SerreCL}, IV, Proposition 14).  Using Lemma \ref{Lbase}, Theorem
\ref{Tsetup}, and induction, it follows that there is a $\ZZ/p^{n-1}$-cover
$Y_{n-1} \to X$ satisfying Claim \ref{claim1} for $L_{n-1}/k[[t]]$, that is \'{e}tale outside $D(r_{n-1})$.  Then
Theorem \ref{Tsetup} (ii) shows that Claim \ref{claim1}, thus Conjecture \ref{OortConj}, holds for $L_n/k[[t]]$. 
\Endproof

\section{The base case} \label{plan2}

The proof of Lemma \ref{Lbase} is straightforward, using the explicit equations given
in \cite{GM98}.

{\bf Proof of Lemma \ref{Lbase} for $\ZZ/p$-extensions:}
By Artin-Schreier theory, any $\ZZ/p$-extension of $k((t))$ is given by an equation $y^p-y=f_1$, where $f_1 \in k((t))$ is unique
up to adding elements of the form $a^p-a$, with $a \in k((t))$.
Thus, we may assume that $f_1=\sum_{i=1}^{m_1}  a_it^{-i}\in k[t^{-1}]$ is a polynomial in
$t^{-1}$ such that $a_i=0$ for $i\equiv 0\pmod{p}$.  Then the break 
in the ramification filtration is $m_1$, which is prime to $p$.  Since 
$\sqrt[m_1]{1/f_1}$ is a uniformizer of $k((t))$, we may make a change of variables and assume $f_1 = t^{-m_1}$.
So we assume our equation is given by $y^p - y = t^{-m_1}$.
In \cite{GM98}, II, Theorem 4.1, a $\ZZ/p$-cover $Y \to \PP^1$ satisfying Claim \ref{claim1} is given by the Kummer extension
\begin{equation}\label{Zplift}
      y^p = G_1(T) = 1 + \lambda^p T^{-m_1},
\end{equation}
where $T$ reduces to $t$ and $\lambda = \zeta_p - 1$ for $\zeta_p$ a primitive $p$th root of unity.
The zeroes of $G_1$ all have valuation $\frac{p}{m_1(p-1)} > \frac{1}{m_1(p-1)} =: r_1$, and the unique pole is at $T=0$.  Since the 
branch points all have valuation greater than $r_1$, the lemma is proved.
\Endproof

{\bf Proof of Lemma \ref{Lbase} for $\ZZ/p^2$-extensions:}
Let $L_2/k((t))$ be a $\ZZ/p^2$-extension with upper ramification breaks $(m_1, m_2),$ and let $r_i = \frac{1}{m_i(p-1)}$ for
$i \in \{1, 2\}$.  After a possible change of variables, \cite{GM98}, II, Theorem 5.5 gives a $\ZZ/p^2$-cover 
$Y \to \PP^1$ satisfying Claim \ref{claim1}, in the form of a Kummer extension 
$$z^{p^2} =  G_1(T)G_2(T)^p.$$  Here $G_1$ is as in \eqref{Zplift} and $G_2$ is a polynomial in $T^{-1}$, which is called
\begin{equation}\label{gmeq}
G(T^{-1}) + p\mu^p \sum_{s=1}^{m_1(p-1)} A_sT^{-s}
\end{equation}
in \cite{GM98}, II, Theorem 5.5.  Also, $v(\mu) = \frac{1}{p(p-1)}$ and $v(A_s) \geq 0$ for all $s$.  
It is clear from the expression for $G(T^{-1})$ given in {\em loc.cit.} that the coefficient of each non-constant term of $G(T^{-1})$ has
valuation at least $\frac{1}{p-1}$.  The same holds for the $p\mu^p A_s$ by inspection.  Furthermore, the proof of
\cite{GM98}, II, Theorem 5.5 shows that each term in (\ref{gmeq}) has degree less than $m_2$ in $T^{-1}$ ($m_2$ is called
$d$ in {\em loc.cit.}).

Now, we have already seen that the zeroes and poles of $G_1$ have valuation greater than $r_1 > r_2$.  The only pole
of $G_2$ is at $T=0$, so it suffices to show that the zeroes of $G_2$ have valuation greater than $r_2$, or equivalently 
(by the theory of Newton polygons), that the coefficient of $T^{-\ell}$ in $G_2 - 1$ has valuation greater than $\ell r_2 = \frac{\ell}{m_2(p-1)}$.  
But this is true because $\ell < m_2$ and the valuation is at least $\frac{1}{p-1}$.
\Endproof

\section{Characters and Swan conductors} \label{Sgeom}

In this section we introduce the general geometric setup (characters) and our
most important tools (Swan conductors). As laid out in the introduction, our
goal is to construct $p^n$-cyclic covers $Y\to X=\PP^1_K$ of the projective
line reducing to a given cover $\Yb\to\Xb=\PP^1_k$ which is \'etale outside
the origin. For technical reasons it is more convenient to work with the
corresponding character of the Galois group of the function field of $X$.

\subsection{Geometric setup}\label{Sgeomsetup}

Let $X$ be a smooth, projective and absolutely irreducible curve over $K$. We
write $\KK$ for the function field of $X$. We assume that $X$ has good
reduction, and fix a smooth $R$-model $X_R$. We let $\Xb:=X_R\otimes_Rk$
denote the special fiber of $X_R$.  We also fix a $K$-rational point $x_0$ on
$X$ and write $\xb_0\in\Xb$ for the specialization of $x_0$ with respect to the
model $X_R$. In our main example, we have $X=\PP^1_K$, $X_R=\PP^1_R$ and
$x_0=0$, but we will not assume this in \S\ref{Sgeom}. 

We let $X^{\rm an}$ denote the rigid analytic space associated to $X$. 
The {\em residue class} of $x_0$ with respect to the model $X_R$,
\[
        D := ]\xb_0[_{X_R} \subset X^{\rm an},
\]
is the set of points of $X^{\rm an}$ specializing to $\xb_0\in\Xb$
(\cite{BLstable}). It is an open subspace of $X^{\rm an}$, isomorphic to the open
unit disk. To make this isomorphism explicit we choose an element
$T\in\OO_{X_R,\xb_0}$ with $T(x_0)=0$ and whose restriction to the special
fiber generates the maximal ideal of $\OO_{\Xb,\xb_0}$ (this is possible
because $\Xb$ is smooth). Then $\hat{\OO}_{X_R,\xb_0}=R[[T]]$, and $T$ induces
an isomorphism of rigid analytic spaces
\[
     D \cong \{\, x\in(\AA^1_K)^{\rm an} \mid v(x)>0 \,\}
\]
which sends the point $x_0$ to the origin. We call $T$ a {\em parameter} for
the open disk $D$ with {\em center} $x_0$. The choice of $T$ having been made,
we identify $D$ with the above subspace of $(\AA^1_K)^{\rm an}$. 

For $r\in\QQ_{\geq 0}$ we define
\[
     D[r] := \{\, x \in D \mid v(x)\geq r \,\}.
\]
We have $D[0]=D$. For $r>0$ the subset $D[r]\subset D$ is an affinoid
subdomain. Let $v_r:\KK^\times\to\QQ$ denote the ``Gauss
valuation" with respect to $D[r]$. This is a discrete valuation on $\KK$ which
extends the valuation $v$ on $K$ and has the property $v_r(T)=r$. It
corresponds to the supremum norm on the open subset $D[r]\subset X^{\rm an}$.

Let $\kappa_r$ denote the residue field of $\KK$ with respect to the valuation
$v_r$. For $r=0$, $\kappa_0$ is the function field of $\Xb$.  We let
$\ord_{\xb_0}:\kappa_0^\times\to\ZZ$ denote the normalized valuation
corresponding to the point $\xb_0\in\Xb$. 

Now suppose that $r>0$. Then after replacing $K$ by a finite extension (which
depends on $r$!) we may assume that $p^r\in K$. Then $D[r]$ is isomorphic to a
closed unit disk over $K$ with parameter $T_r:=p^{-r}T$. Moreover, the residue
field $\kappa_r$ is the function field of the canonical reduction $\Db[r]$ of
the affinoid $D[r]$. In fact, $\Db[r]$ is isomorphic to the affine line over
$k$ with function field $\kappa_r=k(t)$, where $t$ is the image of $T_r$ in
$\kappa_r$.  For a closed point $\xb \in \Db[r]$, we let $\ord_{\xb}:
\kappa_r^\times\to\ZZ$ denote the normalized discrete valuation corresponding
to the specialization of $\xb$ on $\Db[r]$. We let $\ord_\infty$ denote the
unique normalized discrete valuation on $\kappa_r$ corresponding to the `point
at infinity'.

For $F\in \KK^\times$ and $r\in\QQ_{\geq 0}$, we let $[F]_r$ denote the image
of $p^{-v_r(F)}F$ in the residue field $\kappa_r$. 

\subsection{Characters}\label{Scharacters}

We fix $n\geq 1$ and assume that $K$ contains a primitive $p^n$th root of
unity $\zeta_{p^n}$ (this is true after a finite extension of $K$).  For an
arbitrary field $L$, we set
\[
      H^1_{p^n}(L):=H^1(L,\ZZ/p^n\ZZ).
\]
In the case of $\KK$, we have         
\[
      H^1_{p^n}(\KK):=H^1(\KK,\ZZ/p^n\ZZ)
         \cong \KK^\times/(\KK^\times)^{p^n}
\]
(the latter isomorphism depends on the choice of $\zeta_{p^n}$). Elements of
$H^1_{p^n}(\KK)$ are called {\em characters} on $X$.  Given an element $F\in
\KK^\times$, we let $\K_n(F)\in H^1_{p^n}(\KK)$ denote the character
corresponding to the class of $F$ in $\KK^\times/(\KK^\times)^{p^n}$.

For $i=1,\ldots,n$ the
homomorphism
\[
     \ZZ/p^i\ZZ \to \ZZ/p^n\ZZ, \qquad a \mapsto p^{n-i}a,
\]
induces an injective homomorphism $H^1_{p^i}(\KK)\inj H^1_{p^n}(\KK)$. Its
image consists of all characters killed by $p^i$. We consider $H^1_{p^i}(\KK)$
as a subgroup of $H^1_{p^n}(\KK)$ via this embedding.

A character $\chi\in H^1_{p^n}(\KK)$ gives rise to a (possibly branched)
Galois cover $Y\to X$. If $\chi=\K_n(F)$ for some $F\in\KK^\times$, then $Y$
is a connected component of the smooth projective curve given generically by
the Kummer equation $y^{p^n}=F$. If $\chi$ has order $p^i$ as element of
$H^1_{p^n}(\KK)$, then the Galois group of $Y\to X$ is the unique subgroup of
$\ZZ/p^n\ZZ$ of order $p^i$.

A point $x\in X$ is called a {\em branch point} for the character $\chi\in
H^1_{p^n}(\KK)$ if it is a branch point for the cover $Y\to X$. The {\em
  branching index} of $x$ is the order of the inertia group for some point
$y\in Y$ above $x$. The set of all branch points is called the {\em branch
  locus} of $\chi$ and is denoted by $\BB(\chi)$.

\begin{defn}
  A character $\chi\in H^1_{p^n}(\KK)$ is called {\em admissible} if its branch
  locus $\BB(\chi)$ is contained in the open disk $D$.
\end{defn}

\subsubsection{Reduction of characters}

Let $\chi\in H^1_{p^n}(\KK)$ be an admissible character of order $p^n$, and let
$Y\to X$ be the corresponding cyclic Galois cover. Let $G\cong\ZZ/p^n\ZZ$
denote the Galois group of $Y\to X$. Let $Y_R$ be the normalization of
$X_R$ in $Y$. Then $Y_R$ is a normal $R$-model of $Y$ and we have
$X_R=Y_R/G$. 

After enlarging our ground field $K$, we may assume that the character $\chi$
is {\em weakly unramified} with respect to the valuation $v_0$, see
\cite{Epp73}. By definition, this means that for all extensions $w$ of $v_0$
to the function field of $Y$ the ramification index $e(w/v_0)$ is equal to
$1$. It then follows that the special fiber $\Yb:=Y_R\otimes_R k$ is reduced
(see e.g.\ \cite{ArzdorfWewers}, \S 2.2).

\begin{defn}
  We say that the character $\chi$ has {\em \'etale reduction} if the map
  $\Yb\to\Xb$ is generically \'etale. 
\end{defn}

In terms of Galois cohomology the definition can be rephrased as follows. The
character $\chi$ has \'etale reduction if and only if the restriction of
$\chi$ to the completion $\KKh_0$ of $\KK$ with respect to $v_0$ is {\em
  unramified}. The latter means that $\chi|_{\KKh_0}$ lies in the image of the
cospecialization morphism
\[
      H^1_{p^n}(\kappa_0)\to H^1_{p^n}(\KKh_0)
\]
(which is simply the restriction morphism induced by the projection
$\Gal_{\KKh_0}\to\Gal_{\kappa_0}$).  Since the cospecialization morphism is
injective, there exists a unique character $\chib\in H^1_{p^n}(\kappa_0)$
whose image in $H^1_{p^n}(\KKh_0)$ is $\chi|_{\KKh_0}$.  By construction, the
Galois cover of $\Xb$ corresponding to $\chib$ is isomorphic to an irreducible
component of the normalization of $\Yb$.

\begin{defn}
  If $\chi$ has \'{e}tale reduction, we call $\chib$ the {\em reduction} of $\chi$, and $\chi$ a {\em lift} of
  $\chib$. 
\end{defn}

\begin{rem}
  Assume that $\chi$ has \'etale reduction. Then the condition that $\chi$ is
  admissible implies that the cover $\Yb\to\Xb$ corresponding to the reduction
  $\chib$ is \'etale over $\Xb-\{\xb_0\}$ (the proof uses Purity of Branch
  Locus, see e.g.\ \cite{Nagata}). It follows that $\Yb$ is smooth outside the inverse
  image of $\xb_0$.
\end{rem}

\begin{defn}
  Let $\chi\in H^1_{p^n}(\KK)$ be an admissible character of order $p^n$. We
  say that $\chi$ has {\em good reduction} if it has \'etale reduction and the
  cover $\Yb\to\Xb$ corresponding to the reduction $\chib$ of $\chi$ is
  smooth. 
\end{defn}

Note that a $\ZZ/p^n$-cover of $\PP^1_k$, unramified outside $x_0$, is uniquely determined by its germ above the branch point 
(see, e.g., \cite{Ka86}).  Thus, with the above notation, the local Oort conjecture (more specifically, Claim \ref{claim1}) may be reformulated as
follows.
 
\begin{conj} \label{oortconj} Suppose that $X=\PP^1_K$. Let $\chib\in
  H^1_{p^n}(\kappa_0)$ be a character of order $p^n$, unramified outside of
  $\xb_0$.  Then (after replacing $K$ by a finite extension, if necessary)
  there exists an admissible character $\chi\in H^1_{p^n}(\KK)$ with good
  reduction lifting $\chib$.
\end{conj}

\subsection{Swan conductors}\label{swan}

\subsubsection{ } 

Fix $r\in\QQ_{\geq 0}$. We assume that $p^r\in K$. Let $\KKh_r$ denote the
completion of $\KK$ with respect to the valuation $v_r$. Let $\chi\in
H^1_{p^n}(\KK)$ be a character of order $\leq p^n$.  By Epp's theorem
(\cite{Epp73}) we may assume that the restriction $\chi|_{\KKh_r}$ is weakly
unramified.  Under this condition, we can define three types of invariants which
measure in some way the ramification of $\chi$ with respect to the valuation
$v_r$.

First of all, we have the {\em depth Swan conductor}
\[
    \delta_{\chi}(r):=\sw(\chi|_{\KKh_r})\in\QQ_{\geq 0},
\]
see \cite{cyclic}, Definition 3.3.  By definition, $\delta_\chi(r)=0$ if and
only if $\chi$ is unramified with respect to $v_r$. If this is the case then
the reduction $\chib_r\in H^1_{p^n}(\kappa_r)$ is well defined (see the
previous subsection on the case $r=0$).

Let us now assume that $\delta_\chi(r)>0$. Then we can define the {\em
  differential Swan conductor} of $\chi$ with respect to $v_r$,
\[
    \omega_\chi(r):=\dsw(\chi|_{\KKh_r})\in\Omega_{\kappa_r}^1,
\]
see \cite{cyclic}, Definition 3.9. 

Finally, let $\ord_{\xb}:\kappa_r^\times\to\ZZ$ be a normalized discrete
valuation whose restriction to $k$ is trivial. Of course, $\ord_{\xb}$
corresponds either to a closed point $\xb$ on the canonical reduction of the
affinoid $D[r]$, or it corresponds to the point at infinity,
$\xb=\infty$. Then the composite of $v_r$ with $\ord_{\xb}$ is a
valuation on $\KK$ of rank two, which we denote by
$\eta(r,\xb):\KK^\times\to\QQ\times\ZZ$ (see e.g.\ \cite{ZariskiSamuelII}, \S
10, p.\ 43; the group $\QQ\times\ZZ$ is equipped with the lexicographic
ordering). By definition, we have
\[
        \eta(r,\xb)(F) = (v_r(F),\ord_{\xb}([F]_r)),
\]
for $F\in\KK^\times$. In \cite{KatoVC} Kato defines a Swan conductor
$\sw_\chi^{\rm K}(r,\xb)\in\QQ_{\geq 0}\times\ZZ$ of $\chi$ with respect to
$\eta(r,\xb)$ (see \cite{KatoVC}, Definition 2.4 and 3.10; note that we have
$\epsilon:=(0,1)$). By definition, the first component of $\sw_\chi^{\rm
  K}(r,\xb)$ is equal to $\delta_\chi(r)$. We define the {\em boundary Swan
  conductor}
\[
     \sw_\chi(r,\xb) \in\ZZ
\]
as the second component of $\sw_\chi^{\rm K}(r,\xb)$.

\begin{rem} \label{swanrem1}
  The invariant $\sw_\chi(r,\xb)$ is determined by the invariants
  $\delta_\chi(r)$ and $\omega_\chi(r)$, as follows.
  \begin{enumerate}
  \item
    If $\delta_\chi(r)=0$ then
    \[
         \sw_\chi(r,\xb) = \sw_{\chib_r}(\xb).
    \]
    Here $\chib_r$ is the reduction of $\chi$ with respect to $v_r$ (well
    defined because $\chi$ is unramified at $v_r$) and $\sw_{\chib_r}(\xb)$ is
    the usual Swan conductor of $\chib_r$ with respect to the valuation
    $\ord_{\xb}$ (one less than the Artin conductor for nontrivial characters, see \cite{SerreCL}, VI, \S2). 
    This formula follows easily from
    the definitions. As a consequence we see that $\sw_\chi(r,\xb)\geq 0$ and
    that $\sw_\chi(r,\xb)=0$ if and only if $\chib_r$ is unramified with
    respect to $\ord_{\xb}$.
  \item
    If $\delta_\chi(r)>0$ then we have
    \[
         \sw_\chi(r,\xb) = -\ord_{\xb}(\omega_\chi(r))-1.
    \]
    This follows from \cite{KatoVC}, Corollary 4.6.\footnote{The proof in
      \cite{KatoVC} uses class field theory and works only if the residue
      field of $K$ is quasi-finite. For a much more direct and elementary
      proof, see \cite{arizona}.}
  \end{enumerate}
\end{rem}

\begin{prop}\label{Paddchar}
  Let $\chi_1, \chi_2 \in H^1_{p^n}(\KK)$, and let $\chi_3 = \chi_1\chi_2$.
  For $i \in \{1, 2, 3\}$ and $r \in \QQ_{\geq 0}$, set $\delta_i =
  \delta_{\chi_i}(r)$. If $\delta_i>0$ then we set
  $\omega_i:=\omega_{\chi_i}(r)$. If $\delta_i=0$ then $\chib_i\in
  H^1_{p^n}(\kappa_r)$ denotes the reduction of $\chi_i$ with respect to
  $v_r$.  
  \begin{enumerate}
  \item If $\delta_1 \ne \delta_2$ then $\delta_3 = \max(\delta_1, \delta_2)$.
    If $\delta_1 > \delta_2$ then $\omega_3 = \omega_1$.
  \item Assume $\delta_1 = \delta_2>0$. Then
    \[
       \omega_1 + \omega_2 \ne 0 \quad\Rightarrow\quad
          \delta_1=\delta_2=\delta_3,\;\;\omega_3=\omega_1+\omega_2
    \]
    and
    \[
        \omega_1 + \omega_2 = 0 \quad\Rightarrow\quad
           \delta_3<\delta_1.
    \]
  \item Assume $\delta_1 = \delta_2=0$. Then $\delta_3=0$ and
    $\chib_3=\chib_1\chib_2$. 
 \end{enumerate}
\end{prop}

\proof
Parts (i)-(ii) follow from \cite{cyclic}, Proposition 3.10. Part (iii) is
clear, because the cospecialization map $H^1_{p^n}(\kappa_r)\to
H^1_{p^n}(\KKh_r)$ is a homomorphism.
\Endproof

\subsubsection{} 

The finite extension of $K$ that was necessary in order to define the
invariants $\delta_\chi(r)$, $\omega_\chi(r)$ and $\sw_\chi(r,\xb)$ depends on
$r$. However, the values $\delta_\chi(r)$ and $\omega_\chi(r)$ do not depend
on the choice of this extension. Therefore, it makes sense to consider
$\delta_\chi$, $\omega_\chi$ and $\sw_\chi$ as functions in $r\in\QQ_{\geq 0}$
and $\xb$.

\begin{prop}\label{Pdeltalin}
  $\delta_\chi$ extends to a continuous, piecewise linear function 
  \[
      \delta_\chi:\RR\to\RR_{\geq 0 }
  \]
  Furthermore: 
  \begin{enumerate}
  \item
    For $r\in\QQ_{>0}$, the left (resp.\ right) derivative of $\delta_{\chi}$ at $r$ is $-\sw_{\chi}(r, \infty)$ (resp.\ $\sw_{\chi}(r, 0)$).
  \item
    If $r$ is a {\em kink} of $\delta_\chi$ (meaning that the left and right derivatives do not agree), then $r\in\QQ$.
  \end{enumerate}
\end{prop}

\proof See e.g.\ \cite{boundary}, Proposition 2.9. A more direct proof of a
special case of the proposition can be derived from \S\ref{orderp} below.
\Endproof

\begin{cor}\label{Cdeltalin}
If $r \geq 0$ and $\delta_{\chi}(r) > 0$, then the left and right derivatives of $\delta_{\chi}$ at $r$ are given by
$\ord_{\infty}(\omega_{\chi}(r)) + 1$ and $-\ord_0(\omega_{\chi}(r)) - 1$, respectively.
\end{cor}

\proof Immediate from Proposition \ref{Pdeltalin} (i) and Remark \ref{swanrem1} (ii).
\Endproof

\subsubsection{} 

We are going to characterize the case when $\chi$ has good reduction in terms
of the function $\delta_\chi$. Our main tool for this is a certain ``local
vanishing cycles formula". As a special case, we recover the criterion for good
reduction from \cite{GM98}, \S 3.4. We fix an admissible character $\chi\in
H^1_{p^n}(\KK)$ of order $p^n$ and let $Y\to X$ denote the corresponding
Galois cover. Let us also fix $r\in\QQ_{\geq 0}$ and assume that $p^r\in K$.

Suppose first that $r>0$. Then the affinoid subdomain $D[r]\subset X^{\rm an}$
gives rise to an admissible blowup $X_R'\to X_R$ with the following properties
(\cite{ArzdorfWewers} \S3.5 and \cite{BLrigid}): Firstly, $X_R'$ is a semistable curve whose special fiber
$\Xb':=X_R'\otimes_R k$ consists of two smooth irreducible components which
meet in exactly one point. The first component is the strict transform of
$\Xb$, which we may identify with $\Xb$. The second component is the
exceptional divisor $\Zb$ of the blowup $X_R'\to X_R$, which is isomorphic to
the projective line over $k$ and intersects $\Xb$ in the distinguished point
$\xb_0$. By construction, the complement $\Zb^\circ:=\Zb\backslash\{\xb_0\}$ is
identified with the canonical reduction $\Db[r]$ of the affinoid
$D[r]$. In particular this means that the discrete valuation on $\KK$
corresponding to the prime divisor $\Zb\subset X_R'$ is equivalent to the
valuation $v_r$ and that the residue field $\kappa_r$ may be identified with
the function field of $\Zb$.

Let $Y_R'$ denote the normalization of $X_R'$ in $Y$. We obtain a commutative
diagram 
\[\begin{CD}
    Y_R' @>>> Y_R   \\
    @VVV      @VVV  \\
    X_R' @>>> X_R   \\
\end{CD}\]
in which the vertical maps are finite $G$-covers and each horizontal map is the
composition of an admissible blowup with a normalization. Let $\Wb\subset Y_R'$ 
be the exceptional divisor of $Y_R'\to
Y_R$. After enlarging the ground field $K$ we may assume that $\Wb$ is
reduced. Note that this holds if and only if the character $\chi$ is weakly
unramified with respect to the valuation $v_r$, and that this is exactly the
condition we need to define $\delta_\chi(r)$, $\omega_\chi(r)$ and
$\sw_\chi(r,\xb)$. We now choose a closed point $\xb\in \Zb^\circ=\Db[r]$ 
and a point $\yb\in \Wb$ lying over $\Zb$. We let 
\[
      U(r,\xb):=]\xb[_{D[r]}
\]
denote the residue class of $\xb$ on the affinoid $D[r]$. Clearly, $U(r,\xb)$
is isomorphic to the open unit disk. Finally, we let $q:\tilde{\Wb}\to \Wb$ denote
the normalization of $\Wb$ and set 
\[
      \delta_{\yb}:=\dim_k\,(q_*\OO_{\tilde{\Wb}}/\OO_{\Wb})_{\yb}. 
\]
Then $\delta_{\yb}\geq 0$ and we have $\delta_{\yb}=0$ if and only if $\yb\in \Wb$
is a smooth point.

The above notation extends to the case $r=0$ as follows. If $r=0$ then we let
$\Zb:=\Xb$ denote the special fiber of the smooth model $X_R$ of $X$ and
$\Wb:=\Yb$ the special fiber of $Y$. We set $\xb:=\xb_0$ and choose an arbitrary
point $\yb\in \Wb$ above $\xb_0$. The residue class $U(r,\xb)$ is now equal to
the open disk $D$ and the invariant $\delta_{\yb}$ is defined in the same way as
for $r>0$. 

\begin{prop} \label{vcprop}
  With the notation introduced above we have
  \[
        \sw_\chi(r,\xb) = \abs{\BB(\chi)\cap U(r,\xb)}-1-2\delta_{\yb}.
  \]
\end{prop}

\proof This follows from \cite{KatoVC}, Theorem 6.7. To see this, note that
the left hand side of the formula in {\em loc.cit.} (the ``vanishing cycles")
remains invariant if the sheaf $\mathcal{F}$ is pulled back to a Galois cover
of $\Spec(A)$ on which $\mathcal{F}$ becomes constant. In our situation we
take for $A:=\OO_{X_R',\xb}^{\rm h}$ the henselian local ring of $\xb$ on the
scheme $X_R'$ and for $\mathcal{F}$ the \'etale sheaf corresponding to the
character $\chi$. Then the ring extension $B:=\OO_{Y_R',\yb}^{\rm h}/A$ gives
rise to a Galois cover which trivializes $\mathcal{F}$. To prove Proposition
\ref{vcprop}, one applies the formula from {\em loc.cit.} to $\mathcal{F}$ and
its pullback. After equating the left hand side of both formulas and tracing
back the definitions of all terms appearing in the right hand side one obtains
the desired result.
\Endproof

As a first consequence of the above proposition we reprove the following
important criterion for good reduction from \cite{GM98}, \S 3.4.

\begin{cor}  \label{vccor1}
\begin{enumerate}
  \item Let $\chi\in H^1_{p^n}(\KK)$ be an admissible character of order $p^n$. Then
  \[
      \abs{\BB(\chi)} \geq \sw_\chi(0,\xb_0)+1.
  \]
  Also, $\chi$ has good reduction if and only if $\delta_\chi(0)=0$ and
  equality holds above.
\item Suppose $\chi$ has good reduction with upper ramification breaks $(m_1,
  \ldots, m_n)$.  If $1 \leq i \leq n$, then
  \[	
    \abs{\{ x \in\BB(\chi) \, | \, \text{ramification index of } x \text{ is
      exactly } p^{n-i+1}\}} = m_i - m_{i-1},
  \]
  where we set $m_0 = -1$.
\end{enumerate}
\end{cor}

\proof The inequality in part (i) follows immediately from Proposition
\ref{vcprop} since $\BB(\chi)\subset D=U(0,\xb_0)$ by assumption.  Now, by
definition $\chi$ has good reduction if and only if $\delta_\chi(0)=0$ and
$\Yb=W$ is smooth in any point $\yb$ above the distinguished point
$\xb_0$. The latter condition is equivalent to $\delta_{\yb}=0$.  Thus, the
rest of part (i) also follows from Proposition \ref{vcprop}.

In the situation of part (ii), the character $\chi_i := \chi_n^{p^{n-i}} \in
H^1_{p^i}(\KK)$ is an admissible character with reduction $\chib_i$ of order
$p^i$ with upper ramification breaks $(m_1, \ldots, m_i)$.  Thus the Swan
conductor of $\chib_i$ is $m_i$ (\cite{SerreCL}, Corollary 2 to Theorem 1, noting that the Swan conductor is one less than the Artin conductor).  
By part (i), $\abs{\BB(\chi_i)} = m_i + 1$.  Since elements of $\BB(\chi_i)$ correspond
exactly to elements of $\BB(\chi)$ with ramification index at least
$p^{n-i+1}$, part (ii) follows.
\Endproof

\begin{cor} \label{vccor2}
  Let $\chi\in H^1_{p^n}(\KK)$ be an admissible character of order $p^n$,
  $r\in\QQ_{>0}$ and $\xb$ a point on the canonical reduction of $D[r]$. 
  Then
  \[
      \sw_\chi(r,\xb) \leq \abs{\BB(\chi)\cap U(r,\xb)} -1.
  \]
  Moreover, if $\chi$ has good reduction then equality holds. 
\end{cor}

\proof
The inequality follows immediately from Proposition \ref{vcprop}. To prove the
second statement we note that if $\chi$ has good reduction, then in
  the situation of Proposition \ref{vcprop} the point $\yb\in \Wb$ is smooth.
This is because the curve $\Wb$ is the exceptional divisor of the modification $Y_R'\to
Y_R$. If $\chi$ has good reduction then $Y_R$ is smooth over $R$ and hence
regular. It follows from Castelnuovo's criterion (see e.g.\ \cite{LiuAG},
Theorem 9.3.8) that $\Wb$ is smooth. 
\Endproof  

\begin{cor}\label{vccor3}
In the situation of Corollary \ref{vccor2}, if $\delta_{\chi}(r) > 0$, we have $\ord_{\xb}(\omega_{\chi}(r)) \geq -\abs{\BB(\chi)\cap U(r,\xb)}$, with
equality if $\chi$ has good reduction.
\end{cor}

\proof Immediate from Corollary \ref{vccor2} and Remark \ref{swanrem1} (ii).
\Endproof
 
\begin{rem} \label{deltaconcaverem} 
  \begin{enumerate}
  \item If $\chi$ has good reduction, then Corollaries \ref{Cdeltalin} and 
    \ref{vccor3} show that $\delta_\chi$ is a
    piecewise linear, weakly concave down function. Moreover, the position of the kinks of
    $\delta_\chi$ correspond to the valuations of the ramification points in $\BB(\chi).$ If
    $r>0$ is a kink, then the number of ramification points of $\chi$ with
    valuation $r$ is precisely the difference between the left and the right
    derivative of $\delta_\chi$ at $r$.
  \item Now assume that $\BB(\chi)\subset D(r_0)$ for some
    $r_0\in\QQ_{>0}$. Then it follows from Remark \ref{swanrem1},
    Corollary \ref{Cdeltalin}, and Corollary \ref{vccor3} that the
    restriction of $\delta_\chi$ to the interval $[0,r_0]$ is weakly concave up. Together
    with (i) this shows that, if $\chi$ has good reduction, then
    $\delta_\chi|_{[0,r_0]}$ is linear.
  \end{enumerate}
\end{rem}

\subsection{Characters of order $p$} \label{orderp}

\subsubsection{}

We will now describe in the special case $n=1$ how to determine the function
$\delta_\chi$ explicitly in terms of a suitable element $F\in\KK^\times$
corresponding to the character $\chi\in
H^1_p(\KK)\cong\KK^\times/(\KK^\times)^p$. 

\begin{prop} \label{deltaprop} 
  Let $F\in \KK^\times\backslash(\KK^\times)^p$,
  $\chi:=\K_1(F)\in H^1_p(\KK)$ and $r\in\QQ_{\geq 0}$. Suppose that
  $v_r(F)=0$. Suppose, moreover, that $\chi$ is weakly unramified with respect
  to $v_r$. 
  \begin{enumerate}
  \item
    We have 
    \[
        \delta_\chi(r) = \frac{p}{p-1} - \max_H\, v_r(F-H^p),
    \]
    where $H$ ranges over all elements of $\KK$. 
  \item
    The maximum of $v_r(F-H^p)$ in (i) is achieved if and only if
    \[
        g:=[F-H^p]_r\not\in\kappa_r^p.
    \]
    If this is the case, and $\delta_{\chi}(r) > 0$, then
    \[
       \omega_\chi(r) = \begin{cases}
          \;\;dg/g & \text{if $\delta_\chi(r)=p/(p-1)$,} \\
          \;\;dg   & \text{if $0<\delta_\chi(r)<p/(p-1)$.}
                        \end{cases}
    \]
    If, instead, $\delta_{\chi}(r) = 0$, then $\chib$ corresponds to the
    Artin-Schreier extension given by the equation $y^p - y = g$.
  \end{enumerate}
\end{prop}

\proof
The assumption $g\not\in\kappa_r^p$ shows that $v_r(H)\geq 0$. If $v_r(H)>0$
then $g=[F]_r$. In this case, \cite{cyclic}, Proposition 5.3, says that
$\delta_\chi(r)=p/(p-1)$ and $\omega_\chi(r)=dg/g$. Otherwise, 
\[
     FH^{-p} = 1 + p^sG,
\]
where $s:=v_r(F-H^p)>0$ and $G\in\KK$ is an element with $v_r(G)=0$ and
residue class $g$. By \cite{cyclic}, Proposition 5.3, we now have
$\delta_\chi(r)=p/(p-1)-s$. Moreover, if $\delta_\chi(r)>0$ then
$\omega_\chi(r)=dg$.  

If $\delta_{\chi}(r) = 0$, then $s = \frac{p}{p-1}$.  Make a change of
variable $1 + \lambda Y = FH^{-p}$, where $\lambda \in K$ is the unique
element satisfying $\lambda^{p-1} = -p$ and $v(1 +
\frac{\lambda}{p^{1/(p-1)}})>0$.  Then the equation $(1 + \lambda Y)^p = 1 +
p^sG$ yields $Y^p - Y = G + o(1)$, which reduces to the desired Artin-Schreier
extension.
\Endproof

\subsection{Detecting the slope of $\delta_{\chi}$}\label{Skink}

Let $\chi\in H^1_p(\KK)$ be an admissible character of order $p$, giving rise
to a branched cover $Y \to X$.   Let $m > 1$ be a prime-to-$p$ integer.  
We assume that the following conditions hold.
\begin{itemize}
\item[(a)]
  The branch locus of $\chi$ is contained in the closed disk $D[r_0]$, for
  some $r_0>0$, and $T=0$ is one of the branch points.
\item[(b)]
  For all $r \in (0, r_0]$, the left derivative of $\delta_{\chi}$ at $r$ is $\leq m$ (equivalently, by 
  Proposition \ref{Pdeltalin}, $\sw_{\chi}(r, \infty) \geq -m$). 
\item[(c)] For all $r \in (0, r_0]$, we have $\delta_{\chi}(r) > 0$.
\end{itemize}

Because of Condition (a) we can represent $\chi$ as the Kummer class of a
power series
\[
      F=1+\sum_{i=1}^\infty a_iT^{-i},
\]
with $a_i\in R$ and $v(a_i)\geq r_0i$. We wish
to find a {\em polynomial} $H$ in $T^{-1}$ whose $p$th power approximates $F$
well enough to use Proposition \ref{deltaprop} simultaneously for all $r$ in an interval $(0, s] \cap \QQ$, for some $0 < s < r_0$. We will then get explicit expressions for the slopes of
$\delta_\chi$ on the interval $[0,r_0]$.   

For any $N \geq 1$, set
\[
      H:=1+\sum_{j=1}^N b_jT^{-j}.
\]
Here we consider the $b_j$ for the moment as indeterminates. Write
\[
     F-H^p=\sum_{k=1}^\infty c_k T^{-k},
\]
where $c_k$ is a polynomial in $b_1,\ldots,b_{\min(k,N)}$. Note that
$c_k=a_k\in R$ for $k>pN$. 

\begin{lem} \label{kinklem1}
   Assuming condition (a), after replacing $K$ by some finite extension, there exist
  $b_1,\ldots,b_N\in R$ such that
  \begin{enumerate}
  \item
    $v(c_k)\geq r_0k$ for all $k$, and 
  \item
      $c_{kp}=0$ for all $k\leq N$.
  \end{enumerate}
\end{lem}

\proof 
For (ii) to hold we can solve the equations $c_{pN}=c_{p(N-1)}=\ldots=c_p=0$
inductively:
\[\begin{split}
    c_{pN}&=b_N^p-a_{Np}=0,\\
    c_{p(N-1)} &= b_{N-1}^p+\ldots=0,\\
    \vdots &\qquad \vdots
\end{split}\]
One easily checks that these solutions verify $v(b_j)\geq jr_0$ for all
$j$. So we have $v_{r_0}(F),v_{r_0}(H)\geq 0$. Thus
we get $v_{r_0}(F-H^p)\geq 0$, which is equivalent to (i).
\Endproof

\begin{rem}\label{Rkink}
The proof above shows that there are only finitely many solutions for the $b_j$, and that they vary analytically as the 
$a_i$ do.
\end{rem}

\begin{prop}\label{kinkprop} 
Assume conditions (a), (b), and (c) hold.  Choose $s \in (0, r_0) \cap \QQ$ and $N \in \NN$ such that
\begin{equation}\label{kinkeq3}
       pN\geq \frac{p}{(p-1)(r_0-s)}.
\end{equation}
 Let $b_1,\ldots,b_N$ be as in Lemma \ref{kinklem1}. 
 Define $\lambda_m(\chi)\in[0,r_0]$ by
    \[
         \lambda_m(\chi):=\max\big(\{r\in (0,r_0]\mid \sw_\chi(r,\infty)>-m\}
                                \cup\{0\}\big).
    \]
  Set 
  \[
     \mu_m(\chi):=\max\big(\{\frac{v(c_m)-v(c_k)}{m-k} \mid 1\leq k<m\}
                          \cup\{0\}\big).
  \]
  Then 
  \begin{enumerate}
  \item For all $r\in(0,s] \cap \QQ$ we have
  \[
     [F-H^p]_r\not\in\kappa_r^p.
  \]
  Therefore,
  \[
        \delta_\chi(r)= p/(p-1)-v_r(F-H^p), \quad 
          \sw_\chi(r,\infty)= -\ord_\infty[F-H^p]_r.
  \]\
  \item  
    We have
    \[
       \lambda_m(\chi)<s \quad\Leftrightarrow\quad \mu_m(\chi)<s. 
    \]
  \item
    If $\lambda_m(\chi)<s$ then $\lambda_m(\chi)=\mu_m(\chi)$. 
  \end{enumerate}
\end{prop}

\begin{rem}
Note that, if $\lambda_m(\chi) \neq r_0$, then Proposition \ref{Pdeltalin} implies that $\lambda_m(\chi)$ is
the largest value in $(0, r_0]$ where $\delta_{\chi}$ has a kink. 
\end{rem}

\proof (of Proposition \ref{kinkprop}) Fix $r\in(0,s] \cap \QQ$ and set $M:=\ord_\infty[F-H^p]_r$. By definition and
by Lemma \ref{kinklem1} (i) we have
\begin{equation} \label{kinklemeq2}
   v_r(F-H^p)=v(c_M)- rM\geq M(r_0-r)\geq M(r_0-s).
\end{equation}
On the other hand, condition (c) shows that $\delta_\chi(r) >0$ and
Proposition \ref{deltaprop} show that
\begin{equation} \label{kinklemeq3}
  v_r(F-H^p) < p/(p-1).
\end{equation}
Using \eqref{kinklemeq2}, \eqref{kinklemeq3} and the choice of $N$ we obtain
the inequality  
\begin{equation} \label{kinklemeq4}
     M<\frac{p}{(p-1)(r_0-s)}\leq Np.
\end{equation}
If $M$ was divisible by $p$ then \eqref{kinklemeq4} and Lemma \ref{kinklem1}
(ii) would show that $c_M=0$, which contradicts the definition of
$M$. Therefore, $M$ is prime to $p$, and Part (i) of the lemma follows from
Proposition \ref{deltaprop} and Remark \ref{swanrem1} (ii).

In order to prove (ii) and (iii) we note that, by condition (b), 
$\lambda_m(\chi)<s$ is equivalent to $\sw_\chi(s,\infty)=-m$.
Let $\mathcal{N}$ be the Newton polygon of the power series
$F-H^p=\sum_{k=1}^\infty c_kT^{-k}$.  By (i), the Swan conductor of $\chi$ on $(0,s]$ is
determined by $\mathcal{N}$. In particular, we have $\sw_\chi(s,\infty)=-m$ for
some $r\in(0,s]$ if and only if the point $(m,v(c_m))$ is a vertex of
$\mathcal{N}$ and $s_1\leq s < s_2$, where $s_1$ (resp.\ $s_2$) is the
slope of the edge to the left (resp.\ to the right) of the vertex $(m,v(c_m))$.
Furthermore, in this case we have $\mu_m(\chi) = s_1 < s$, and 
$\lambda_m(\chi) = s_1$, which proves (iii).

To prove (ii), it remains to show that if $\mu_m(\chi) < s$, then $\lambda_m(\chi) < s$.  
If $(m, v(c_m))$ is a vertex of $\mathcal{N}$, then as in the paragraph above 
we have $\lambda_m(\chi) = \mu_m(\chi),$ so $\lambda_m(\chi) < s$.
If $(m, v(c_m))$ is not a vertex of $\mathcal{N}$, then there exists a line segment of
$\mathcal{N}$ with slope $s' < \mu_m(\chi) < s$ connecting two points $(i, v(c_i))$ and $(j, v(c_j))$,
where $i < m < j$.  But this means that $\sw_{\chi}(s, \infty) \geq -j$, which contradicts 
condition (b).  This proves (ii).
\Endproof

Proposition \ref{kinkprop} will be the key to \S\ref{A3}.

\section{Proof of Theorem \ref{Tsetup}}\label{Sproof}
\subsection{Plan of the proof} \label{plan}

We continue with the notation of \S\ref{Sgeom}, and for the rest of the paper, we set $X \cong \PP^1$ and $x_0 = 0$.  
Recall that $D$ is the unit disk in $(\AA^1_K)^{\rm an}$ centered at $0$.  For $r \in \QQ_{\geq 0}$, we set
$$D(r) = \{T \in (\AA^1_K)^{\rm an}\, | \, \abs{T} < \abs{p}^r \} \subseteq D.$$
We are given a character $\chib_n \in
H^1_{p^n}(\kappa_0)$ of order exactly $p^n$, with upper ramification breaks $(m_1, m_2, \ldots, m_n)$.  We further assume that
$n \geq 2$.  For $1 \leq i \leq n$, set $r_i = \frac{1}{m_i(p-1)}$.  Recall that $p \nmid m_1$ and
\[
       m_i\geq p\,m_{i-1},
\]
for $i=2,\ldots,n$. Moreover, if the inequality above is strict then
$(m_i,p)=1$.  For $i=1,\ldots,n$ we set $\chib_i:=\chib_n^{p^{n-i}}\in H^1_{p^i}(\kappa_0)$. 
By assumption, for each $1 \leq i < n$, there is a character $\chi_{i}$ lifting 
$\chib_{i}$.  We assume that $\BB(\chi_{n-1})$ lies in the disk $D(r_{n-1})$, and we may further assume that 
$T=0$ is a branch point of order $p^{n-1}$.  In order to prove Theorem \ref{Tsetup}, we must 
show that there exists a character $\chi_n \in H^1_{p^n}(\KK)$ with (good) reduction 
$\chib_n$.  Furthermore, we must have $\BB(\chi_n) \subseteq D$, and if there is no integer $a$ satisfying 
\begin{equation}\label{aeq}
\frac{m_n}{p} - m_{n-1} < a \leq \left(\frac{m_n}{m_n - m_{n-1}}\right) \left(\frac{m_n}{p} - m_{n-1}\right),
\end{equation}
then we must even have $\BB(\chi_n) \subseteq D(r_n)$. 
We will construct $\chi_n$ such that $\chi_n^p = \chi_{n-1}$.

We may assume that $\chi_{n-1}$ corresponds to an extension of $\KK$ 
given by a system of Kummer equations
\[
   y_i^p=y_{i-1}G_i, \quad i=1,\ldots,n-1
\]
with $y_0:=1$ and $G_i \in \KK$.  Any $\chi \in H^1_{p^n}(\KK)$ such that $\chi^p = \chi_{n-1}$ is given by an additional equation
\begin{equation}\label{gn}
  y_n^p=y_{n-1}G.
\end{equation}
Since we must have $\BB(\chi) \subseteq D$,
we will search for $G \in 1+ T^{-1}\m[T^{-1}],$ where $\m$ is the maximal ideal of $R$.  
In particular, 

\begin{equation}\label{Egnform}
     G=\prod_{i=1}^N(1-x_iT^{-1})^{a_i},
\end{equation}
where $a_i\in\NN$, $(a_i,p)=1$, and $x_i\in \m$ are pairwise
distinct. We will say that the polynomial $G$ \emph{gives rise}
to the character $\chi$.  If $x_i$ is a branch point of $\chi_{n-1}$
then we may also transfer the term $(1-x_it^{-1})^{a_i}$ into
$G_{n-1}$. Therefore, we may assume that none of the $x_i$ is a
branch point of $\chi_{n-1}$.  If this is the case, then Corollary \ref{vccor1} (ii) shows that a necessary condition
for good reduction of $\chi$ is that $N=m_n-m_{n-1}$.  We assume this.  

Our proof that a choice $G_n$ for $G$ exists giving rise to a character $\chi_n$ whose (good) reduction is $\chib_n$ will 
be done in two parts: 
\begin{itemize}
\item[{(\bf Part A)}] 
  We prove that there exists a polynomial $G_{\min} \in 1 + T^{-1}\m[T^{-1}]$ 
  giving rise to a character $\chi_{\min}$ 
  with good reduction $\chib_{\min}$ having upper ramification breaks $(m_1, \ldots, m_{n-1}, pm_{n-1})$ at the ramified point.
\item[{(\bf Part B)}]
  We construct the desired polynomial $G_n$ by modifying $G_{\min}$. 
\end{itemize}
Furthermore, we show that if there is no $a$ satisfying \eqref{aeq}, and if $G_n$ gives rise to $\chi_n$, then
$\BB(\chi_n) \subseteq D(r_n)$.  This will complete the proof of Theorem \ref{Tsetup}.

We remark that the basic strategy for our proof is adapted from the proof of
the case $n=2$ by Green and Matignon, \cite{GM98}. Essentially, Part (A)
corresponds to Lemmas 5.2 and 5.3 in {\em loc.cit.}, whereas Part (B) corresponds to
Lemma 5.4.   

The proof of Part (A) will be done in three steps.  The first step (\S\ref{A1}) is to find an appropriate family of candidate polynomials for $G_{\min}$,
which we will call $\G_n$.  This family is defined in Definition \ref{Adef1}.  
The second step (\S\ref{A2}) is to show that if a polynomial in $\G_n$ yields a character with bad reduction, it can be 
altered (within $\G_n$) to obtain a new polynomial whose reduction is ``closer" to being good (i.e., the depth Swan conductor of the
corresponding character is lower).  The key result here is Proposition \ref{Aprop1}.   
Lastly (\S\ref{A3}), we show there must exist a polynomial in $\G_n$ that is ``closest" to having good reduction (i.e., the depth Swan conductor of the
corresponding character is minimal).  This is the content of Proposition \ref{kinkminprop}.  
Combining these steps shows that there must exist $G_{\min} \in \G_n$ giving a character with
good reduction.

Proposition \ref{Pmain} proves Part (B), and is found in \S\ref{PartB}.

\subsection{A family of candidate polynomials} \label{A1}
\subsubsection{} \label{A1.1}

We continue with the setup of \S\ref{plan}.  In particular, recall that $\chi_i$ is a lift of
$\chib_i$ for $1 \leq i < n$, and $\chi$ is the character arising from $G$, as in \eqref{gn}.  
If $r \in \QQ_{\geq 0}$, then to simplify the notation we will write $\delta_i(r)$ and
$\omega_i(r)$ instead of $\delta_{\chi_i}(r)$ and $\omega_{\chi_i}(r)$ for the depth and differential Swan conductors
of $\chi_i$ (\S\ref{swan}).  Furthermore,
write $\delta_n(r)$ and $\omega_n(r)$ instead of $\delta_{\chi}(r)$ and $\omega_{\chi}(r)$.   
As will become apparent in Proposition \ref{prop1} and its proof, it will be very important to control 
$\omega_n(r_{n-1})$.  The polynomials $G$ which give our desired $\omega_n(r_{n-1})$ will comprise our
candidate family $\G_n$.

\begin{lem} \label{swanlem}
  Assume $1 \leq i<n$. Then for $r\in [0,r_i]$ we have
  \begin{equation} \label{swanlem1eq1}
       \delta_i(r) = m_i\cdot r.
  \end{equation}
  Moreover, for $0<r\leq r_i$ we have 
  \begin{equation} \label{swanlem1eq2}
        \omega_i(r) = \frac{c_i\,dt}{t^{m_i+1}}.
  \end{equation}
  Here $c_i\in k^\times$ is a constant depending on $i$ and $r$. 
\end{lem}

\proof By hypothesis, $\chi_i$ has good reduction $\chib_i$. Therefore, $\delta_i(0)=0$. 
On the other hand, the hypothesis that
$\chi_i$ has good reduction and that all of its $m_i+1$ branch points are
contained in the disk $D(r_i)$ implies that
\[
      \ord_0(\omega_i(r)) = -m_i-1, \quad \ord_{\xb}(\omega_i(r)) = 0 
\]
for all $r\in(0,r_i]$ and $\xb \ne 0, \infty$, using Corollary \ref{vccor3}.  So \eqref{swanlem1eq1} follows, using
Corollary \ref{Cdeltalin}. But now the same corollary shows that
\[
     \ord_\infty(\omega_i(r)) = m_i-1,
\]
and \eqref{swanlem1eq2} follows as well.
\Endproof

\begin{rem} \label{swanrem}
  Suppose $1<i<n$ and $0<r<r_i$. Lemma \ref{swanlem} shows that 
  \[
     p\delta_{i-1}(r)\leq \delta_i(r)< 1/(p-1).
  \]
  Moreover, the first inequality is an equality if and only if
  $m_i=pm_{i-1}$. It follows from \cite{cyclic}, Theorem 4.3 (ii) that 
  $\C(\omega_i(r))=0$ if and only if $m_i>pm_{i-1}$, where $\C$ is the Cartier operator. This is consistent with
  \eqref{swanlem1eq2}. On the other hand, if
  $m_i=pm_{i-1}$ then 
  \[
        \C(\omega_i(r)) = \omega_{i-1}(r).
  \]
  In particular, we have $c_i=c_{i-1}^p$ in \eqref{swanlem1eq2}. 
\end{rem}

\subsubsection{} \label{A1.4}

We now focus on the critical radius $r_{n-1}$.  To further simplify the notation
we will, until the end of \S\ref{A1}, write $\omega_i$ (resp.\ $\delta_i$) instead of
$\omega_i(r_{n-1})$ (resp.\ $\delta_i(r_{n-1})$).  By Lemma
\ref{swanlem} we have $\delta_{n-1}=1/(p-1)$. So \cite{cyclic}, Theorem 4.3,
says that
\begin{equation} \label{Aeq4}
  \delta_n = \frac{p}{p-1}
\end{equation}
and
\begin{equation} \label{Aeq5}
  \C(\omega_n) = \omega_n+\omega_{n-1}.
\end{equation} 
Let $m$ be the minimal upper ramification break $m_i$ such that $m_{n-1}$ is a power of $p$ times $m_i$.  Thus $m$ is prime to $p$.
Set $\nu = n-1-i$.  Thus $0 \leq \nu \leq n-2$, and $m_{n-1} = mp^{\nu}$.

By Lemma \ref{swanlem} and Remark \ref{swanrem} we have
\[
   \omega_i=\frac{c\,dt}{t^{m+1}},\ldots,
      \omega_{n-1} = \frac{c^{p^{\nu}}dt}{t^{p^{\nu}m+1}},
\]
for some $c \in k^{\times}$.
After a change of parameter we may assume that $c=m$, viewed as an element of $k^{\times}$.  Note that $m^p = m$.  Set
\begin{equation}\label{eta}
    \eta:=-(\omega_i+\ldots+\omega_{n-1})
         = -m\sum_{j=0}^{\nu}t^{-mp^j-1}dt.
\end{equation}

\begin{lem} \label{Alem1} 
  Let $g = [G]_{r_{n-1}}$.  Then
  \[
      \omega_n = \eta +\frac{dg}{g}.
  \]
\end{lem}

\proof 
One easily checks that $\C(\eta)=\eta+\omega_{n-1}$. Using \eqref{Aeq5}
we conclude that
\[
     \omega_n=\eta+\frac{dh}{h},
\]
for some $h\in\kappa_{r_{n-1}}^\times$. 

Let us first assume that $G=1$. Then $\BB(\chi_n)$ lies in
the disk $D(r_{n-1})$. It follows from Corollary \ref{vccor3} that the differential
$\omega_n$ has no poles outside $t=0$. Since $\eta$ has no poles outside
$t=0$, this can happen only if $dh/h=0$. This shows that the lemma is true if
$G=1$.

To prove the general case, we note that multiplying $G$ by an element
$H\in\KK^\times$ has the effect of adding the character $\K_1(H)\in
H^1_p(\KK)$ to $\chi$. We may assume that $v_{r_{n-1}}(H)=0$ and let
$h\in\kappa_{r_{n-1}}$ denote the residue of $H$. Then we have
\[
    \delta_{\K_1(H)} = p/(p-1) \quad\text{and}\quad
    \omega_{\K_1(H)} = dh/h
\]
if and only if $h\not\in\kappa_{r_{n-1}}^p$, by Proposition
\ref{deltaprop}. If $h\in\kappa_{r_{n-1}}^p$ then
$\delta_{\K_1(H)}<p/(p-1)$. In both cases, Proposition \ref{Paddchar} shows
that multiplying $G$ by $H$ has the effect of adding $dh/h$ to $\omega_n$. The lemma follows.
\Endproof

The following proposition is not strictly necessary for the proof of Theorem \ref{Tsetup}, but it helps to narrow our search for the correct
$G_n$.  Recall that we assume $G_n$ to be in the form \eqref{Egnform}, and that $N = m_n - m_{n-1}$.
\begin{prop} \label{prop1}
  If $\chi$ has good reduction then the following hold.
  \begin{enumerate}
  \item
    For all $i$ we have $v(x_i)\leq r_{n-1}= \frac{1}{m_{n-1}(p-1)}$.
  \item
    For $i,j$ with $v(x_i)=v(x_j)=r_{n-1}$ we have $\bar{x}_i\neq\bar{x}_j$
    (where $\bar{x}_i$ denotes the reduction of $x_ip^{-r_{n-1}}$).
  \item
    Write $N=N_1+N_2$, where $N_1$ is the number of $x_i$'s with
    $v(x_i)=r_{n-1}$. We may assume that $v(x_i)<r_{n-1}$ for
    $i=N_1+1,\ldots,n$. Then 
    \[
        \sum_{i=N_1+1}^N\,a_i \equiv 0 \pmod{p}.
    \]
    In particular, if $a_i=1$ for all $i$ then $N_2\equiv 0\pmod{p}$. 
  \item
    If $m_n=pm_{n-1}$ then $N_1=m_{n-1}(p-1)$ and $N_2=0$. Otherwise,
    $N_1<m_{n-1}(p-1)$ and $N_2>0$.
  \end{enumerate}
\end{prop}

\proof
By Lemma \ref{Alem1}, we have
\begin{equation} \label{eq1.3}
   \omega_n = dg/g-m\sum_{j=0}^{\nu}t^{-mp^j-1}\,dt,
\end{equation}
where
\[
    g:=[G]_{r_{n-1}}\in k[t^{-1}].
\]
It follows that $\ord_0(\omega_n)=-m_{n-1}-1$. So if $\chi$ has good
reduction, then Corollary \ref{vccor3} shows that the number of branch points specializing to $0$ (i.e.\ with
valuation $>r_{n-1}$) must be equal to $m_{n-1}+1$. Since $\chi_{n-1}$ has
exactly $m_{n-1}+1$ branch points with valuation $>r_{n-1}$, none of the new
branch points can have this property. This proves (i).

By \eqref{eq1.3}, $\omega_n$ can have at most a simple pole at any 
point $\bar{x}\neq 0$, and then good reduction and Corollary \ref{vccor3} implies that branch points with radius $r_{n-1}$ have to
lie in distinct residue classes. This proves (ii).  It follows similarly from Corollary \ref{vccor3} that $\omega_n$ has no zeroes outside
$t = \infty$, so the fact that $\text{div}(\omega_n)$ has degree $-2$ means that $\ord_\infty(\omega_n)=m_{n-1}+N_1-1 \geq 0$. 
But it is easy to see that 
\[
    \ord_\infty(g)=\sum_{i=N_1+1}^N\,a_i,
\]
and \eqref{eq1.3} shows that $\ord_\infty(\omega_n) \geq -1$ if
and only if $\ord_{\infty}g \equiv 0\pmod{p}$.  This proves (iii). 

On the other hand we have $\C(\omega_n)=\omega_n+\omega_{n-1}$ and
$\ord_\infty(\omega_{n-1})=m_{n-1}-1$, which implies, by an easy calculation,
that 
\[
    \ord_\infty\omega_n\leq pm_{n-1}-1.
\]
It follows that $N_1\leq m_{n-1}(p-1)$. 

Now suppose that $m_n=pm_{n-1}$.  Then Corollaries \ref{Cdeltalin} and \ref{vccor3} show that
the right derivative of $\delta_n$ is at most $m_{n-1} + N = pm_{n-1}$ on $[0, r_{n-1})$.
Since $\delta_n(r_{n-1}) = \frac{p}{p-1}$ by \eqref{Aeq4}, 
and good reduction requires $\delta_n(0) = 0$, this slope must be $pm_{n-1}$ on the entire interval.  Thus $\ord_\infty\omega_n
= m_n -1 = pm_{n-1}-1$, by Corollary \ref{Cdeltalin}.  Hence $N_1=N=m_{n-1}(p-1)$ and $N_2=0$. Otherwise, if
$m_n>pm_{n-1}$, then the condition that $\delta_n$ is weakly concave down (Remark \ref{deltaconcaverem}) and has right derivative $m_n$
at $r=0$ (Proposition \ref{Pdeltalin} (i)) implies that $\delta_n(r)>pm_{n-1}r$ for $0 < r<r_{n-1}$. But this means
that $\ord_\infty\omega_n+1<pm_{n-1}$ at $r=r_{n-1}$, hence
$N_1<m_{n-1}(p-1)$. It follows that
\[
    N_2=m_n-m_{n-1}-N_1>m_n-pm_{n-1}>0.
\]
This completes the proof of the proposition.
\Endproof

It follows from Proposition \ref{prop1} that, up to a constant factor that we may eliminate by rescaling $t$, we have
\begin{equation}\label{Egred}
     [G]_{r_{n-1}} = g=t^{a_0}\prod_{i=1}^{N_1} (1-\bar{x}_it^{-1})^{a_i},
\end{equation}
where 
\[
    a_0:=-\sum_{i=N_1+1}^N a_i \equiv 0 \pmod{p}.
\]
Hence
\begin{equation}\label{Elaurent}
   dg/g = \sum_{i=1}^{N_1} \frac{a_i\bar{x}_it^{-2}dt}{1-\bar{x}_it^{-1}}.
\end{equation}

\begin{cor}\label{Cgform}
  In the notation of Proposition \ref{prop1}, if $\chi$ has good reduction then
  \begin{equation} \label{eq1.7}
      \omega_n = dg/g-m\sum_{j=0}^{\nu}t^{-mp^j-1}dt = 
          \frac{c\,dt}{t^{m_{n-1}+1}\prod_{i=1}^{N_1}(t-\bar{x}_i)},
  \end{equation}
  where $c:=(-1)^{N_1+1}m (\prod_{i=1}^{N_1}\bar{x}_i$) is a nonzero constant.
  In particular, $\ord_{\infty} \omega_n = m_{n-1} + N_1 - 1.$
\end{cor}

\proof The middle expression is the expression deduced for $\omega_n$ in
\eqref{eq1.3}. We have seen in the proof of Proposition \ref{prop1} that $\omega_n$
has simple poles at the $\bar{x}_i$, a pole of order $m_{n-1}+1$ at $0$, and no
zero outside $\infty$. It follows that $\omega_n$ is equal to the right hand
side of \eqref{eq1.7} times a constant. To determine this constant, one computes the Laurent
series representation in $t$ of both sides. 
\Endproof

The next theorem, showing that we can often find a $g$ satisfying the conditions of \ref{Cgform}, is critical.  

\begin{thm} \label{Athm1} 
  Suppose $m | m_n$ (equivalently, $m | N = m_n - m_{n-1}$).  
  Then, under the assumption that $a_i = 1$ for $i \geq 1$, there is a solution $g \in k[t^{-1}]$ to \eqref{eq1.7} as in 
  \eqref{Egred} of degree $N$.
  \end{thm}

\proof This is contained in Corollary \ref{1tom}.
\Endproof 

\begin{defn} \label{Adef1} 
  If $g$ is the solution to \eqref{eq1.7} guaranteed by Theorem \ref{Athm1}, then we define $\G_n$ to be the subset of all $G$ of the
  form \eqref{Egnform} with $[G]_{r_{n-1}}=g$. 
\end{defn}

\begin{rem}\label{Rn1}
Corollary \ref{1tom} also shows that, if $a_i  = 1$ for all $i \geq 1$ and $N_1$ is as in \eqref{Egred}, then we have $m | N_1$ and 
$$m_{n-1}(p-1)-mp<N_1\leq m_{n-1}(p-1), \quad N_1\equiv N\pmod{p}.$$  This determines $N_1$ uniquely.
This means that if $G$ is as in \eqref{Egnform}, has all $a_i = 1$, and gives rise to $\chi$ with good reduction, then the number of 
zeroes of $G$ with valuation $r_{n-1}$ is fixed.  The importance of this condition is illustrated in Example \ref{exa2}.
\end{rem}

\begin{exa} \label{exa2} 
  Assume that $p=5$, $n=3$, and the upper ramification breaks of $\chib_3$ are $(m_1, m_2, m_3) = (1, 5, 34)$.   
  If $G_3$ is in the form of \eqref{Egnform} with all $a_i = 1$, and if $G_3$ gives rise to a character $\chi_3$ with good reduction $\chib_3$, 
  then $N = 29$.  Remark \ref{Rn1} shows that $\chi_3$ gives rise to a cover with exactly $N_1=19$
  branch points at radius $r_2=1/20$ and $N_2=10$ branch points at radius
  $<1/20$.   We know from \eqref{Aeq4} that $\delta_3(1/20) = \frac{5}{4}$.  By Corollaries \ref{Cdeltalin} and \ref{vccor3},
  the (left) slope of $\delta_3$ at $r \leq 1/20$ is equal to $5 + N(r)$, where
  $N(r)$ is the number of branch points with valuation $\geq r$.  Since $\delta(0) = 0$, one can show that there must be at
  least one branch point with valuation $\leq 1/200$. But $r_3=1/136>1/200$, so $\BB(\chi_3)$ cannot lie in the disk 
  $D(r_3)$.  
  
  See Remark \ref{needineq} and Example \ref{Econditionaldisk} for such an example with
  $10$ branch points with valuation exactly $1/200$.  Note that this does not contradict Theorem \ref{Tsetup} (ii), as we can take
  $a = 2$ there.
\end{exa}

\begin{rem}\label{Rarb}
If we do not assume $a_i = 1$ for all $i$, then \eqref{eq1.7} can still sometimes be solved.  
In particular, in light of Example \ref{exa2}, it would be nice to find solutions to \eqref{eq1.7} with arbitrary $a_i$ and $N_1$  
not satisfying Remark \ref{Rn1}.  We might then have some hope of finding a lift $\chi_n$ of $\chib_n$ with $\BB(\chi_n)$ lying in the disk
$D(r_n)$, even when the condition in Theorem \ref{Tsetup} (ii) does not hold.
However, even when such solutions to \eqref{eq1.7} exist, it seems as if our current techniques are often insufficient to turn them into lifts.
For further discussion, see Remark \ref{samebounds}.
\end{rem} 

\subsection{Reducing the depth Swan conductor}\label{A2}
We maintain the notation of \S\ref{A1}, and we assume further that
$m_n = pm_{n-1}$.  Then, by Proposition \ref{prop1}, we have $N = N_1 = m_{n-1}(p-1)$.  

Recall that any $G \in\G_n$ (Definition \ref{Adef1})
gives rise to a character $\chi$ of order $p^n$ lifting $\chi_{n-1}$ as in \eqref{gn}, by adjoining the equation
$y_n^p = y_{n-1}G$.  
By Corollaries \ref{Cgform} and \ref{Cdeltalin}, we know that the left derivative of $\delta_n$ at $r_{n-1}$ is $m_n$.
Recall also from \eqref{Aeq4} that $\delta_n(r_{n-1}) = \frac{p}{p-1} = m_nr_{n-1}$.  
It follows that there exists $0\leq \lambda<r_{n-1}$ such that
\[
     \delta_n(s) = s\cdot m_n = psm_{n-1} = p\delta_{n-1}(s)
\]
for all $s\in[\lambda,r_{n-1}]$ (the last equality following from Lemma \ref{swanlem}). Let $\lambda(G)$ be the minimal value of
$\lambda$ with this property.  In other words, $\lambda(G)$ is the largest kink
of the function $\delta_n$ on the open interval $(0,r_{n-1})$ (or is zero if
$\delta_n$ is linear on $[0,r_{n-1}]$).
Note that $\delta_n(\lambda) = m_n\lambda < \frac{p}{p-1}$.

\begin{prop}\label{Pzerogood}   
If $G \in \G_n$ satisfies $\lambda(G) = 0$, then the corresponding character $\chi$ has good reduction.
\end{prop}

\proof By definition, $\lambda(G) = 0$ implies $\delta_{\chi}(0) = 0$.  Now,
If $\chi$ lifts $\chi_{n-1}$ and has $\delta_{\chi}(0) = 0$, then
$\sw_{\chi}(0, \xb_0) \geq pm_{n-1}$, as the Swan conductor of a $p^n$-cyclic
extension must be at least $p$ times the Swan conductor of its index $p$
subextension (because each upper ramification break of a $p^n$-cyclic
extension must be at least $p$ times the previous one).  On the other hand, by
construction, $\BB(\chi)$ has exactly $N + m_{n-1} +1 = pm_{n-1} +1$ branch
points.  The proposition then follows by Corollary \ref{vccor1} (i).
\Endproof

Thus, in order to prove Part (A), it suffices to show that $\lambda(G)=0$ for some
$G \in\G_n$.  Proposition \ref{Aprop1} will show that $\lambda(G)=0$ is the
only possible minimal value of $\lambda(G)$, as $G$ ranges over $G_n$.  
In \S\ref{A3}, we will show that this minimum is realized.  First, we state a lemma.

\begin{lem} \label{Alem3} Let $G \in\G_n$ and $r\in[0,r_{n-1})$. Let
  $f\in t^{-1}k[t^{-1}]$ be a polynomial of degree $< m_n$ without constant
  term, which we regard as an element of $\kappa_r$. Set
  $s:=p/(p-1)-m_nr$.  Then, after a possible finite extension of $K$, there exists $G'\in\G_n$ and
  $F\in\KK$ such that $v_r(F) = 0$, $[F]_r = f$, and
  \[
       G'/G \equiv 1-p^sF \pmod{(\KK^\times)^p}.
  \]
\end{lem} 

\proof
The proof is given as Corollary \ref{proof}.
\Endproof

\begin{prop} \label{Aprop1} Suppose $G \in\G_n$ with $\lambda(G)>0$. Then there
  exists $G'\in\G_n$ with $\lambda(G')<\lambda(G)$.
\end{prop}

\proof
If $\lambda:=\lambda(G)>0$ then by Corollary \ref{Cdeltalin} we have
that $\ord_\infty(\omega_n(\lambda)) + 1$ is the left derivative of $\delta_n(r)$.  Since $\delta_n$ is concave up at $\lambda$ (Remark \ref
{deltaconcaverem} (ii)), we conclude that $$\ord_{\infty}(\omega_n(\lambda)) < m_n - 1.$$  By hypothesis we have
$p\,\delta_{n-1}(\lambda)=\delta_n(\lambda) = m_n\lambda < p/(p-1)$.  Therefore \cite{cyclic}, Proposition 4.3 (ii) shows that
\[
    \C(\omega_n(\lambda))=\omega_{n-1}(\lambda)=\frac{c\,dt}{t^{m_{n-1}+1}},
\]
for some $c \in k$.  It follows that 
\[
   \omega_n(\lambda)=\frac{c^p\,dt}{t^{m_n+1}} + df,
\]
for some $f\in\kappa_\lambda$. Note that, by Corollaries \ref{Cdeltalin} and \ref{vccor3},
\[
    \ord_0(\omega_n(\lambda))=-m_n-1,\quad \ord_{\xb}(\omega_n(\lambda))\geq 0 \;\;\forall\, \xb\neq 0.
\]
We may therefore assume that $f$ is a polynomial in $t^{-1}$ of
degree $< m_n$ and without constant term.  

By Lemma \ref{Alem3}, there exists $G'\in\G_n$ such that
\[
    G'/G \equiv 1-p^sF \pmod{(\KK^\times)^p},
\]
where $v_{\lambda}(F) = 0$, where $[F]_{\lambda} = f$, and where
$s:=p/(p-1)-m_n\lambda$. Replacing $G$ by the polynomial $G'$ has the effect
of adding $\psi:=\K_1(G'/G)$ to $\chi$. Using Proposition
\ref{deltaprop}, we see that
\[
     \delta_{\psi}(\lambda)=p/(p-1)-s=\delta_n(\lambda).
\]
Therefore, Proposition \ref{Paddchar} shows that the effect on
$\omega_n(\lambda)$ is addition of $-df$ and the result is that
\[
     \ord_{\infty}(\omega_n(\lambda)) = m_n+1.
\]
We conclude, using Corollary \ref{Cdeltalin}, that $\lambda(G')<\lambda(G)$. 
\Endproof

\begin{rem}\label{nobranch}
  An important reason we must assume that $\BB(\chi_{n-1}) \in D(r_{n-1})$ is
  to ensure that no branch point of $\chi_{n-1}$ has valuation less than
  $\lambda$.  If there were such a branch point, then
  $\ord_{\infty}(\omega_n(\lambda))$ above could be negative, which would
  allow $f$ not to be a polynomial in $t^{-1}$, which would prevent us from applying Lemma \ref{Alem3}.
\end{rem}

\subsection{The minimal depth Swan conductor}\label{A3} 
We continue with the notation of \S\ref{A2}, as well as the assumption that
$m_n = pm_{n-1}$ and all $a_i = 1$ for $i \geq 1$.  
To finish the proof of Part (A) from \S\ref{plan}, we must show that the
function $\lambda:\G_n \to \QQ_{\geq 0}$ defined in \S\ref{A2} takes the value $0$ for some $\chi
\in \G_n$.  By Proposition \ref{Aprop1}, the existence of such a $\chi$ is
established by the following proposition.

\begin{prop} \label{kinkminprop} 
  The function $\chi\mapsto \lambda(\chi)$
  takes a minimal value on $\G_n$.
\end{prop}

The rest of \S\ref{A3} is devoted to the proof of this proposition.
\subsubsection{A lemma from rigid analysis} \label{Srigidlem}

The following lemma, which is an easy consequence of the {\em maximum modulus
  principle}, is a crucial ingredient in the proof of Proposition
\ref{kinkminprop}.

\begin{lem} \label{rigidlem}
  Let $X=\Spm(A)$ be an affinoid domain over $K$ and $f_1,\ldots,f_n\in A$
  analytic functions on $X$. Then the function
  \[
       \phi:X\to\RR, \quad \phi(x):=\max_{1\leq i\leq n}\sqrt[i]{\abs{f_i(x)}}
  \]
  takes a minimal value. Equivalently, the function
  \[
       x\mapsto \min_i \frac{v(f_i(x))}{i}
  \]
  takes a maximal value on $X$. 
\end{lem}

\proof
Let $B/A$ be a finite ring extension which contains elements $g_i\in B$ such
that $g_i^i=f_i$, for $i=1,\ldots,n$. Then $B$ is again an affinoid
$K$-algebra, and the induced morphism $q:Y:=\Spm(B)\to X$ is finite and
surjective. For any point $y\in Y$ we have
\[
      \phi(q(y))=\max_{1\leq i\leq n}\abs{g_i(x)}.
\]
So by \cite{BGR}, Lemma 7.3.4/7, the function $\phi\circ q$ takes its minimal
value on $Y$. Since $q$ is surjective, this shows that $\phi$ takes its
minimal value on $X$.
\Endproof

\subsubsection{An affinoid containing $\G_n$}
Let $\G$ be the set of all polynomials of the form
\[
     G= \prod_{i=1}^{N_1}(1-x_iT^{-1}),
\]
where $x_i\in\Kb$ has valuation $r_{n-1}$ and where the residue classes of
$x_i/p^{r_{n-1}}$ in $k$ are all distinct.  Given $G \in \G$, we may assume (after
passing to a finite extension of $K$) that $x_1,\ldots,x_{N_1}\in K$. In this
way, we can consider $\G$ as a subset of affine $N_1$-space over $K$ via the
coordinates $x_i$ and $G$ as a $K$-rational point.
It is clear that $\G$ is an affinoid subdomain of $(\AA^{N_1}_K)^{\rm an}$. 
We identify elements of $\G$ with the characters $\chi$ that they give rise to (\S\ref{plan}). 

Recall that it is a consequence of Proposition \ref{prop1} (iv) that $\G_n \subseteq \G$.  In particular,
\[
       \G_n =\{G \in\G \mid [G]_{r_{n-1}}=g \},
\]
where $g$ is the unique solution of
Equation \eqref{eq1.7} guaranteed by Theorem \ref{Athm1}. As a rigid analytic space,
$\G_n$ is isomorphic to the open unit polydisk.  The idea of the proof is to show that $\lambda(\G)$
takes a minimal value on $\G$, and that the point where this minimum is achieved must lie in $\G_n$.

\subsubsection{}\label{kinkmin1}
     
Let $\phi_{n-1}:Y_{n-1}\to X$ be the Galois cover corresponding to the
character $\chi_{n-1}$. By our induction hypothesis, it has good
reduction and is totally ramified above $T=0$. It follows that the rigid analytic
subspace $C:=\phi_{n-1}^{-1}(D)\subseteq Y_{n-1}$ is an open disk and contains the
unique point $y_{n-1}\in Y_{n-1}$ above $T=0$. We choose a parameter
$\tilde{T}$ for the disk $C$ such that $\tilde{T}(y_{n-1})=0$. Then 
\[
     T=\tilde{T}^{p^{n-1}}u(\tilde{T}), \quad
      \text{with $u(\tilde{T})\in R[[\tilde{T}]]^\times$.}
\]
We conclude that for $r>0$ the inverse image of the closed disk $D(r)\subset
D$ defined by the condition $v(T)\geq r$ is the closed disk $C(\tilde{r})$
defined by $v(\tilde{T})=\tilde{r}:=p^{-n+1}r$. Set
$\tilde{r}_{n-1}:=p^{-n+1}r_{n-1}$. Let $\KK_{n-1}$ denote the function field
of $Y_{n-1}$. 

Let us fix, for the moment, $\chi\in \G$ such that $\chi^p = \chi_{n-1}$. Let
$\tilde{\chi}:=\chi|_{\KK_{n-1}}\in H^1_p(\KK_{n-1})$ denote the restriction
of $\chi$ to $\KK_{n-1}$. If $\chi$ corresponds to a cover $Y\to X$, then
$\tilde{\chi}$ corresponds to the cover $Y\to Y_{n-1}$.  If $\chi \in \G_n$, then in analogy to $\lambda(\chi)$, 
we write $\lambda(\tilde{\chi})$ for the minimum
$\tilde{r} \in [0, \tilde{r}_{n-1}]$ such that $\delta_{\tilde{\chi}}$ is linear on $[\tilde{r}, \tilde{r}_{n-1}]$.
If $\chi \in \G \backslash \G_n$, then we define $\lambda(\tilde{\chi}) = \tilde{r}_{n-1}$.

\begin{lem} \label{kinklem4}\ 
  \begin{enumerate}
  \item For $\chi \in \G_n$, we have $\lambda(\tilde{\chi})=p^{-n+1}\lambda(\chi)$.
  \item  Let 
  $$\tilde{m} = p^n m_{n-1} - \sum_{i=1}^{n-1} m_i(p-1)p^{i-1}.$$ 
  Then $p \nmid \tilde{m}$ and the character
    $\tilde{\chi}\in H^1_p(\KK_{n-1})$ satisfies the conditions (a), (b), and (c) 
    of \S \ref{Skink} (with respect to $\tilde{m}$, the open disk $C\subset
    Y_{n-1}^{\rm an}$ and the family of subdisks $C(\tilde{r})$,
    $\tilde{r}\in[0,\tilde{r}_{n-1}]$).
  \item If $\lambda_{\tilde{m}}(\tilde{\chi})$ is as in Proposition \ref{kinkprop}, and if we set $r_0$ in Proposition \ref{kinkprop} equal to $\tilde{r}_{n-1}$,  
  then $\lambda(\tilde{\chi}) = \lambda_{\tilde{m}}(\tilde{\chi})$ for
  all $\chi \in \G$.
  \end{enumerate}
\end{lem}

\proof  For $r>0$ we
systematically use the notation $\tilde{r}:=p^{-n+1}r$. Then the valuation
$v_{\tilde{r}}$ on $\KK_{n-1}$ (corresponding to the Gauss norm on
$C(\tilde{r})$) is the unique extension of $v_r$.  By \cite{cyclic}, \S7.1 we have 
\[\begin{split}
    \delta_{\tilde{\chi}}(\tilde{r}) & =\psi_{\KK_{n-1}/\KK}(\delta_\chi(r))\\
     & = \delta_\chi(r)-\big(\delta_1(r)\frac{p-1}{p^{n-1}}+\ldots
             +\delta_{n-1}(r)\frac{p-1}{p}\big),
\end{split}\]
where $\psi$ is the inverse Herbrand function (\cite{SerreCL}, IV, \S3). Since all the characters $\chi_i$ ($1 \leq i < n$) have
good reduction and their branch points are contained in $D(r_{n-1})$, it
follows from Remark \ref{deltaconcaverem} (ii) that each $\delta_i$ ($1 \leq i < n$) is linear of slope $m_i$ on the
interval $[0,r_{n-1}]$. Therefore, we have
\[
    \delta_{\tilde{\chi}}(\tilde{r}) = \delta_{\chi}(r) + (\frac{\tilde{m}}{p^{n-1}} - pm_{n-1})r. 
\]
Thus, the left slope of $\delta_{\tilde{\chi}}$ at $\tilde{r}$ is equal to $p^{n-1}c + \tilde{m} - p^nm_{n-1}$, where $c$ is the left-slope of 
$\delta_{\chi}$ at $r$.  Part (i) follows immediately.  Part (ii) follows from the fact that $c \leq m_n = pm_{n-1}$.  
Part (iii) also follows from this fact, along with Proposition \ref{Pdeltalin} and the fact that $\sw_{\chi}(r_{n-1}, \infty) = pm_{n-1}$ iff $\chi \in \G_n$.
\Endproof

Explicitly, the character $\tilde{\chi}$ is the Kummer class of the element
\[
   F:=G_1^{1/p^{n-1}}G_2^{1/p^{n-2}}\cdot\ldots\cdot G_{n-1}^{1/p}G\in\KK_{n-1}^\times.
\]
We write $F$ as a power series in the parameter $\tilde{T}$:
\[
      F=1+\sum_{\ell=1}^\infty a_\ell\tilde{T}^{-\ell}
\]
Note that, since $G_1, \ldots, G_n$ are fixed, $F$ is uniquely determined by the choice of $G$.  So we may consider
the coefficients $a_\ell$ as functions on the space $\G$. It is easy to see
that the $a_\ell$ are analytic functions on $\G$ which are bounded by $1$. In
fact, $a_\ell$ is a polynomial in the coordinates $x_i$ with coefficients in
$R$. 

\subsubsection{} \label{kinkmin2}

We continue with the proof of Proposition \ref{kinkminprop}. Let
$\chi_0\in\G_n$ be an arbitrary lift lying in the residue class
determined by the reduction $g$. From the discussion at the beginning of \S\ref{A2}, it follows that
$\lambda(\chi_0)<r_{n-1}$. We may therefore
choose a rational number $s\in(\lambda(\chi_0),r_{n-1})$. Recall from \S\ref{kinkmin1} that $\tilde{\chi}$ is the restriction of
$\chi$ to the function field $\KK_{n-1}$ of $Y_{n-1}$.  Then by Lemma
\ref{kinklem4} we have
\[
     \lambda(\tilde{\chi})<\tilde{s}:=p^{1-n}s<\tilde{r}_{n-1}.
\] 
We also choose an integer $N$ such that
\[
    Np \geq \frac{p}{(p-1)(\tilde{r}_{n-1}-\tilde{s})}.
\]
Compare with \eqref{kinkeq3}. 

\begin{lem} \label{kinkminlem2}
  There exists a finite cover $\G'\to\G$ and analytic functions
  $b_1,\ldots,b_N$ on $\G'$ with the following property. Set
  \[
       H:=1+\sum_{j=1}^N b_j\tilde{T}_1^{-j}
  \]
  and write
  \[
     F-H^p=\sum_{\ell=1}^\infty c_\ell\tilde{T}_1^{-\ell},
  \]
  where the $c_\ell$ are now analytic functions on $\G'$. Then:
  \begin{enumerate}
  \item
    For all $\ell\geq 1$ and all points $x\in\G'$ we have $v(c_\ell(x))\geq
    r_0\ell$. 
  \item
    We have $c_{p\ell}=0$ for $\ell\leq N$.
  \end{enumerate} 
\end{lem}

\proof
By Lemma \ref{kinklem1} and Remark \ref{Rkink}, there are finitely many solutions for the $b_j$ at each point in $\G$ and the solutions
vary analytically as the $a_{\ell}$ vary in $\G$.  This proves the lemma.
\Endproof

\subsubsection{}

We can now complete the proof of Proposition \ref{kinkminprop}. 
Let $\tilde{m}$ be as in Lemma \ref{kinklem4}, and define the
function $\mu_{\tilde{m}}:\G'\to\RR$ by the formula
\[
    \mu_{\tilde{m}}(x):= \max\big(\{\frac{v(c_{\tilde{m}}(x))-v(c_\ell(x))}{\tilde{m}-\ell}\mid 1\leq \ell <\tilde{m}\}
      \cup\{0\}\big).
\]
Let $\chi\in\G$, write $\tilde{\chi}:=\chi|_{\KK_{n-1}}$ for its restriction
to the function field of $Y_{n-1}$, and let $x\in\G'$ be an arbitrary point above
$\chi$.  By Lemma \ref{kinkminlem2}, we can apply Proposition \ref{kinkprop}
to compare $\mu_{\tilde{m}}(x)$ to $\lambda_{\tilde{m}}(\tilde{\chi})$, which by Lemma \ref{kinklem4} is
equal to $\lambda(\tilde{\chi})$. We conclude that
$\mu_{\tilde{m}}(x)<\tilde{s}$ if and only if $\lambda(\tilde{\chi})<\tilde{s}$. Moreover,
if this is the case then we have $\mu_{\tilde{m}}(x)=\lambda(\tilde{\chi})$. Note also
that in any case we have $\lambda(\chi)=p^{n-1}\lambda(\tilde{\chi})$ when $\chi \in \G_n$.

We apply these arguments twice. Firstly, let $\chi_0\in\G_n$ be the
character with $\lambda(\chi_0)<s$ from the beginning of \S
\ref{kinkmin2}. Let $\chi_0'\in\G'$ be a point above $\chi_0$. Then
$\mu_{\tilde{m}}(\chi_0')<\tilde{s}$.
  
It follows from Lemma \ref{rigidlem} that the function $\mu_m$ takes a minimum
on $\G'$. Let $x\in\G'$ be a point where this minimum is achieved, and let
$\chi\in\G$ be the corresponding lift. We have
$\mu_{\tilde{m}}(x)\leq\mu_{\tilde{m}}(\chi_0')<\tilde{s}$.  Since $$\lambda_{\tilde{m}}(\tilde{\chi}) = \mu_{\tilde{m}}(x) < s < \tilde{r}_{n-1}$$ 
we see that $\chi \in \G_n$.  Applying the above arguments a second time, we
conclude that $\lambda(\chi)=p^{n-1}\mu_{\tilde{m}}(x)$, and that this is actually the
minimal value of the function $\lambda:\G_n \to\RR$.  This completes the proof
of Proposition \ref{kinkminprop}.
\Endproof

Combining Propositions \ref{Pzerogood}, \ref{Aprop1}, and \ref{kinkminprop} finishes the proof of Part (A) from
\S \ref{plan}.  So there is a polynomial
$G_{\min} \in \G_n$ giving rise to a character $\chi_{\min}$ with good reduction and upper ramification breaks $(m_1, \ldots, m_{n-1}, pm_{n-1})$
at the ramification point.

\subsection{Beyond minimality}\label{PartB}

We now prove Part (B) from \S\ref{plan}.  Maintain the notation of the
previous parts of \S\ref{Sproof}.  Let $G_{\min} \in \G_n$ be such that its
corresponding character $\chi_{\min}$ has good reduction $\chib_{\min}$.  Such
a $G_{\min}$ exists by Part (A).  Note that $G_{\min}$ is a polynomial in $1 +
T^{-1}\m[T^{-1}]$ of degree $m_{n-1}(p-1)$.

Recall that $\chib_n$ is our original character, with upper ramification breaks $(m_1, \ldots, m_n)$, and that $m_n$ is \emph{not} necessarily equal to
$pm_{n-1}$.  Furthermore, we saw in \S\ref{Switt} that $\chib_n$ corresponds (upon completion at $t=0$) to a (truncated) Witt vector 
$w_n := (f_1, \ldots, f_n) \in W_n(k((t)))$, and we may assume that each $f_i$ is a polynomial in
$k[t^{-1}]$, all of whose terms have prime-to-$p$ degree.  Then \eqref{rambreaks} shows that $m_n = \max(pm_{n-1}, \deg(f_n))$.  
On the other hand, $\chib_{\min}$ has upper ramification breaks $(m_1, \ldots, m_{n-1}, pm_{n-1})$, and corresponds to a Witt vector 
$w_{\min} := (f_1, \ldots, f_{n-1}, f_{\min})$, where $f_{\min} \in k[t^{-1}]$ has degree $< pm_{n-1}$ and only terms of prime-to-$p$ degree.  
Subtracting Witt vectors yields
\begin{equation}\label{witt}
w_n - w_{\min} = (0, \ldots, 0, f_n - f_{\min}).
\end{equation}
Let $f = f_n - f_{\min}$, which has degree $\leq m_n$, and let $F \in T^{-1}R[T^{-1}]$ be such that $\deg F = \deg f$, 
$v_0(F) = 0$, and $[F]_0 = f$.  

\begin{prop}\label{Pnonminimal}
After a possible finite extension of $K$, there exists $\epsilon \in \QQ_{>0}$, 
as well as $G_n \in 1 + T^{-1}\m[T^{-1}]$ of degree at most $m_n - m_{n-1}$ and 
$H \in 1 + T^{-1}\m[T^{-1}]$ such that $(H, G_n)$ is a solution to 
$$G_{\min}H^p - G_n \equiv -p^{\frac{p}{p-1}}F \pmod{p^{\frac{p}{p-1} + \epsilon}R[T^{-1}]}.$$ 
\end{prop}

\proof This follows immediately from Corollary \ref{Cmodifiedsolution}, substituting $m_n$ and $r_n$ for $m_n'$ and $r_n'$.
\Endproof

\begin{prop}\label{Pwithindisk}
Suppose that there is no integer $a$ satisfying
\begin{equation}\label{noa}
\frac{m_n}{p} - m_{n-1} < a \leq \left(\frac{m_n}{m_n - m_{n-1}}\right)\left(\frac{m_n}{p} - m_{n-1}\right).
\end{equation}
Then we can find $H$ and $G_n$ as in Corollary \ref{Cmodifiedsolution} such that
$v_{r_n}(G_n-1) > 0$, where $r_n = \frac{1}{m_n(p-1)}.$  
Thus all zeroes of $G_n$ lie in the open disk $D(r_n)$.
\end{prop}

\proof This follows immediately from Proposition \ref{Pconditionaldisk}, substituting $m_n$ and $r_n$ for $m_n'$ and $r_n'$.
\Endproof

\begin{prop}\label{Pmain}
There is a character $\chi_n$ whose (good) reduction is $\chib_n$, such that $\BB(\chi_n) \subseteq D$.  If there is no $a \in \ZZ$ satisfying
\eqref{noa}, then $\BB(\chi_n) \subseteq D(r_n)$, where $r_n = \frac{1}{m_n(p-1)}$.  
\end{prop}

\proof
The character $\chi_n$ will correspond to $G_n$, in the notation of Proposition \ref{Pnonminimal}.  From that proposition, we have 
\begin{equation}\label{minton}
\frac{G_n}{G_{\min}}H^{-p} = 1 + p^{\frac{p}{(p-1)}}\frac{F}{G_{\min}H^p} =: \tilde{F}.
\end{equation}
Now, $v_0(\tilde{F} - 1) = \frac{p}{p-1}$ and $[\tilde{F} - 1]_0 = f$, which is not a $p$th power.  Thus, Proposition \ref{deltaprop} shows
that if $\chi_{\tilde{F}} = \K_1(\tilde{F}) \in H_p^1(\KK)$, then $\delta_{\chi_{\tilde{F}}}(0) = 0$ and the reduction $\chib_{\tilde{F}}$ corresponds 
to the Artin-Schreier extension given by $y^p - y = f$.  
Thus if $\chi_{\tilde{F}}' = \K_n(\tilde{F}^{p^{n-1}}) \in H_{p^n}^1(\KK)$, then we also have $\delta_{\chi_{\tilde{F}}'}(0) = 0$ and
the reduction $\chib_{\tilde{F}}'$ corresponds to the same extension, which
is encoded by the Witt vector $(0, \ldots, 0, f)$.

On the other hand, note that $\chi_{\min}$ corresponds to the element 
$$G_1G_2^p\cdots G_{n-1}^{p^{n-2}}G_{\min}^{p^{n-1}} \in \KK^{\times}/(\KK^{\times})^{p^n},$$ 
whereas the character $\chi_n$ coming from $G_n$ corresponds to the element
$$G_1G_2^p\cdots G_{n-1}^{p^{n-2}}G_n^{p^{n-1}} \in \KK^{\times}/(\KK^{\times})^{p^n}.$$
By \eqref{minton}, we have
$\chi_{\min}\chi_{\tilde{F}}' = \chi_n$ as elements of $H^1_{p^n}(\KK)$.  
By Proposition \ref{Paddchar} (iv), we have that $\chib_{\min}\chib_{\tilde{F}}'$ is the reduction of $\chi_n$.
Since $\chib_{\tilde{F}}'$ corresponds to the Witt vector $(0, \ldots, 0, f)$, it follows from \eqref{witt} that the reduction of $\chi_n$ is, in fact, $\chib_n$.
Since $G_n$ is a polynomial of degree $\leq m_n - m_{n-1}$ in $T^{-1}$, we have that $|\BB(\chi_n)| \leq m_n+1$.
By Corollary \ref{vccor1} (i), we have equality, and thus $\chi_n$ has good reduction $\chib_n$, proving the 
first assertion of the proposition.  The second assertion then follows immediately from Proposition \ref{Pwithindisk}.
\Endproof

Proposition \ref{Pmain} completes the proof of Part (B), and Theorem \ref{Tsetup} follows immediately.

\begin{rem}\label{needineq}
Example \ref{Econditionaldisk} shows that it is possible for the result of Proposition \ref{Pwithindisk} not to hold when 
there is an $a \in \ZZ$ satisfying \eqref{noa} (in particular, we take $p=5$, $n=3$, $(m_1, m_2, m_3) = (1, 5, 34)$, and $a=2$).
When this is the case, the branch locus of $\chi_n$ generated above does not lie in the disk $D(r_n)$.
\end{rem}

\section{Proofs of technical results}\label{Stechnical}
In this section, we give the proofs of Theorem \ref{Athm1}, Lemma \ref{Alem3}, and Propositions \ref{Pnonminimal} and \ref{Pwithindisk}.
In fact, we will prove Theorem 
\ref{Athm1} and Lemma \ref{Alem3} in somewhat more generality.  All the proofs are related to each other and will share much notation.  

Throughout \S\ref{Stechnical}, 
we will use notation parallel to that used in \S\ref{Sproof}.  Let $(m_1, \ldots, m_n)$ be a sequence of positive integers such that $p \nmid m_1$, 
that $m_i \geq pm_{i-1}$ for $1 \leq i \leq n$, and that if $p | m_i$, then $m_i = pm_{i-1}$.  For $1 \leq i \leq n$, set $r_i = \frac{1}{m_i(p-1)}$.
Write $N = m_n - m_{n-1}$.  Let $N_1$ and $N_2$ be nonnegative integers
such that $N_1 + N_2 = N$ and $p | N_2$.  Lastly, let $m$ be the minimal $m_i$ such that $m_{n-1}$ is a $p$th power times $m_i$.  
Set $\nu = n-1-i$.  Thus, $m_{n-1} = mp^{\nu}$.

\subsection{Arbitrary types}\label{Sarb}
We work under two helpful assumptions.

\begin{ass} \label{ass1}
  There exist $a_0, a_1,\ldots,a_{N_1}\in\ZZ$ and
  $\bar{x}_1,\ldots,\bar{x}_{N_1}\in k^\times$ with $p|a_0$, with $0 < a_i <p$ for $i \geq 1$, and with $\bar{x}_i\neq\bar{x}_j$,
  such that the differential form
    \[
       \omega:=dg/g-m\sum_{j=0}^{\nu}t^{-mp^j-1}dt, \qquad\text{with}\;\;
            g:= t^{a_0}\prod_{i=1}^{N_1}(1-\bar{x}_it^{-1})^{a_i},
    \]
    satisfies $\ord_\infty(\omega)=m_n-N_2-1 = N_1 + m_{n-1} - 1$ (cf.\ \eqref{eq1.3}).

\end{ass}

To formulate the second assumption we need more notation. We have, by \eqref{Elaurent}, that
\[
    dg/g = \sum_{i=1}^{N_1}\frac{a_i\bar{x}_it^{-2}dt}{1-\bar{x}_it^{-1}}
         = \sum_{\ell=0}^\infty \Big(\sum_{i=1}^{N_1} a_i\bar{x}_i^{\ell+1}\Big) t^{-\ell-2}dt.
\]
Assumption \ref{ass1} is therefore equivalent to the system of
equations 
\[  
\sum_ia_i\bar{x}_i^{\ell+1} =\,
    \begin{cases}
       \;\;m, & \ell=mp^j-1, \ \ 0 \leq j \leq \nu\\ \;\;0, & \text{otherwise,}
    \end{cases}  
\]
for $\ell=0,\ldots,m_n-N_2-2=N_1+m_{n-1}-2$. The Jacobi matrix of this system of
equation at the point $(\bar{x}_i)$ is the $(N_1+m_{n-1}-1,N_1)$-matrix
\[
      \big((\ell+1)a_i\bar{x}_i^\ell\big)_{\ell,i}
\]
over $k$.  The rows of this matrix corresponding to an index $\ell$ with $\ell\equiv -1\pmod{p}$
vanish. Crossing out these trivial rows we obtain the matrix 
\[
      \big((\ell+1)a_i\bar{x}_i^\ell\big)_{\ell\not\equiv -1 (\text{mod } p),i}.
\]
For our key result, Lemma \ref{improvelem}, we need this matrix to be invertible; this is the case if
and only if the matrix
\begin{equation} \label{eq3.6}
  A:= \big(\bar{x}_i^\ell\big)_{\ell\not\equiv -1 (\text{mod } p),i}
\end{equation}
is invertible. 

\begin{ass} \label{ass2}
  The matrix $A$ in \eqref{eq3.6} is invertible.
\end{ass}

\begin{rem}\label{squarerem}
Note that a trivial necessary condition for Assumption \ref{ass2} is that $A$ is a
square matrix.  This is not in general true.  For instance, if $m = 1$, so that $m_{n-1} = p^{n-2}$, then $A$ is square if and only if
$$p^{n-1} - p^{n-2} - p < N_1 \leq p^{n-1} - p^{n-2}.$$ 
\end{rem}

Fix $g$ as in Assumption \ref{ass1}, and assume Assumption \ref{ass2}.  In the field $K(T)$, set $\tilde{T} = p^{-r_{n-1}}T$.  
Let $\G_n$ denote the family of polynomials in $R[T^{-1}]$ of the form
\[
      G=\prod_{i=1}^{N_1} (1-x_i\tilde{T}^{-1})^{a_i},
\]
such that $x_i$ reduces to $\xb_i$.  In particular, for any $G \in \G_n$, we have $[G]_{r_{n-1}} = t^{-a_0}g$.
Lastly, let $m' = N_1 + m_{n-1}$.

\begin{lem} \label{improvelem} 
Under Assumptions \ref{ass1} and \ref{ass2}, let $G\in\G_n$, and let $J \in 1 + \tilde{T}^{-1}\m[\tilde{T}^{-1}]$. 
\begin{enumerate}
\item There exists a unique $G' \in \G_n$ and a unique polynomial $I \in 1 +  \tilde{T}^{-p}\m[\tilde{T}^{-p}]$ of degree $\leq m'-1$ in 
$\tilde{T}^{-1}$
such that $$\frac{G'}{G}I \equiv J \pmod{\tilde{T}^{-m'}}.$$  If $J \equiv 1 \pmod{p^{\sigma}, \tilde{T}^{-m'}}$ for 
$\sigma \in \QQ_{> 0}$, then $v_{r_{n-1}}(\frac{G'}{G}-1) \geq {\sigma}$ and $v_{r_{n-1}}(I-1) \geq {\sigma}$.
\item Let $0 < s < \frac{p}{p-1}$ be a rational number.  After a possible finite extension of $K$,
there exists $G' \in \G_n$ and a polynomial $H \in 1 + \tilde{T}^{-1}\m[\tilde{T}^{-1}]$ of degree $\leq [(m'-1)/p]$ 
such that we have
$$\frac{G'}{G}H^p \equiv J \pmod{p^s, \tilde{T}^{-m'}}.$$
If $J \equiv 1 \pmod{p^{\sigma}, \tilde{T}^{-m'}}$ for some $0 < {\sigma} < \frac{p}{p-1}$, then we can choose $G'$ and $H$ such that 
$v_{r_{n-1}}(\frac{G'}{G} - 1) \geq {\sigma}$ and $v_{r_{n-1}}(H^p - 1) \geq {\sigma}$.
\end{enumerate}
\end{lem}     
     
\proof 
By assumption we have
\[
      G=\prod_{i=1}^{N_1} (1-x_i\tilde{T}^{-1})^{a_i},
\]
and where $x_i\in R$ is a lift of $\bar{x}_i$. We set
\[
    G'=\prod_{i=1}^{N_1} (1-x_i'\tilde{T}^{-1})^{a_i}, \quad x_i':=x_i+z_i,
\]
and where the $z_i$ are for the moment considered as indeterminates. We also set
\[
      I := 1 + \sum_{j=1}^{[(m'-1)/p]} b_j \tilde{T}^{-pj},
\]
for another system of indeterminates $b_j$. Write
\[
      \frac{G'}{G}I = 1 +\sum_{\ell=1}^\infty c_\ell \tilde{T}^{-\ell},
\]
where $c_\ell$ is a formal power series in $(z_i,b_j)$. A simple computation
shows that 
\begin{equation} \label{eq3.8} 
\begin{split}
  \frac{\partial c_\ell}{\partial z_i}|_{z_i=b_j=0} =
        & \;a_ix_i^{\ell-1}  \\
  \frac{\partial c_\ell}{\partial b_j}|_{z_i=b_j=0} = &\;
    \begin{cases} \;\;1, & \ell=pj \\ \;\;0, & \ell\neq pj. \end{cases}
\end{split}
\end{equation} 
The congruence 
\begin{equation} \label{eq3.9}
  \frac{G'}{G}I \equiv J \pmod{\tilde{T}^{-m'}}
\end{equation}
corresponds to a system of $m'-1$ equations (one equation for each monomial
$c_{\ell}\tilde{T}^{-\ell}$, $\ell=1,\ldots,m'-1$) in the indeterminates $(z_i,b_j)$. The Jacobi matrix of this system of
equations is invertible over $R$ if and only if its reduction is invertible over $k$.  From
\eqref{eq3.8} it is easy to see that this is true iff the matrix $A$ from \eqref{eq3.6} is
invertible, which is the case by Assumption \ref{ass2}.  We conclude that
\eqref{eq3.9} has a unique solution with $z_i,b_j\in\m$.  In fact, by the effective Hensel's Lemma, $v(z_i)$ and $v(b_j)$
are all at least as large as $v_{r_{n-1}}(J - 1)$.  This proves (i). 

To prove (ii), we will build $G'$ and $H$ through successive approximation.
Let $G'_1$ and $I_1 = 1 + \sum_{j=1}^{[(m'-1)/p]} b_{j,1} \tilde{T}^{-pj}$ be the unique solution guaranteed by (i).  
So $\frac{G'_1}{G}I_1 \equiv J \pmod{\tilde{T}^{-m'}}$.
Set
\[
      H_1:=1+\sum_{j=1}^{[(m'-1)/p]} b_{j,1}^{1/p}\tilde{T}^{-j}
\]
for any choice of $p$th roots, and set $J_1 := \frac{G_1'}{G}H_1^p$.  Since $H_1^p \equiv I_1 \pmod{p}$, we have that 
$\frac{J}{J_1} \equiv 1 \pmod{p, \tilde{T}^{-m'}}$.  Now, let $G_2'$ and $I_2 = 1 + \sum_{j=1}^{[(m'-1)/p]} b_{j,2} \tilde{T}^{-pj}$ 
be the unique solution to $$\frac{G_2'}{G_1'}I_2 \equiv \frac{J}{J_1} \pmod{\tilde{T}^{-m'}}$$ guaranteed by (i).  
Note that, since the coefficients
of $\frac{J}{J_1} \pmod{\tilde{T}^{-m'}}$ have valuation at least $1$, Part (i) gives that
$v(b_{j,2}) \geq 1$ for all $j$.  Let 
$$H_2 := 1 + \sum_{j=1}^{[(m'-1)/p]} b_{j,2}^{1/p}\tilde{T}^{-j}.$$  Then $H_2^p \equiv I_2 \pmod{p^{1 + 1/p}}$.  Thus
\begin{equation}\label{eqapprox}
\frac{G_2'}{G}(H_1H_2)^p = \frac{G_2'}{G_1'}H_2^p \cdot \frac{G_1'}{G} H_1^p \equiv \frac{J}{J_1}\cdot J_1 \equiv J \pmod{p^{1 + 1/p}, 
\tilde{T}^{-m'}}.
\end{equation} 

We can repeat this process for all $i \in \NN$, letting $G_i'$ and $I_i$ be the unique solution to
$$\frac{G_i'}{G_{i-1}'}I_i \equiv \frac{J}{J_{i-1}} \pmod{\tilde{T}^{-m'}}$$ guaranteed by (i), and constructing $H_i$ from $I_i$ in the same manner 
as $H_1$ and $H_2$.  Set $J_i = \frac{G_i'}{G}(H_1\cdots H_i)^p.$ 
If $\gamma_i = \sum_{j=1}^{i-1}1/p^j$, then $\frac{G_i'}{G_{i-1}'}H_i^p \equiv \frac{J}{J_{i-1}} \pmod{p^{\gamma_i}}.$  
 Analogously to \eqref{eqapprox}, one derives
$$\frac{G_i'}{G}(H_1 \cdots H_i)^p = \frac{G_i'}{G_{i-1}'} H_i^p \cdot \frac{G_{i-1}'}{G}(H_1\cdots H_{i-1})^p 
\equiv J \pmod{p^{\gamma_i}, \tilde{T}^{-m'}}.$$ 
Since $s < \frac{p}{p-1}$, we have that $\gamma_i > s$ for large enough $i$.  Setting $G' = G_i'$ and $H = (H_1 \cdots H_i)$ for
such an $i$ gives the desired solution to part (ii).

To check the last statement of part (ii), it suffices to show that $v_{r_{n-1}}(H_i^p-1) \geq {\sigma}$ for all $i \in \NN$.
Since $v_{r_{n-1}}(I_1-1) \geq {\sigma}$ by (i), we have $v_{r_{n-1}}(H_1^p-1) \geq {\sigma}$ 
(here we use that $\sigma < \frac{p}{p-1}$ to estimate the cross terms).  
Since part (i) also shows $v_{r_{n-1}}(\frac{G_1'}{G}-1) \geq {\sigma}$, we see that 
$v_{r_{n-1}}(J_1-1) \geq {\sigma}$.  An easy induction now shows that $v_{r_{n-1}}(H_i^p-1) \geq \sigma$ for all $i$.
\Endproof

\begin{cor}\label{Capprox}
Under Assumptions \ref{ass1} and \ref{ass2}, 
let $G \in \G_n$, let $r \in \QQ \cap (0, r_{n-1})$, and let $f \in t^{-1}k[t^{-1}]$ be a polynomial of degree $< m' = m_{n-1} + N_1$ 
without constant term.  Fix some $s \in \QQ$ with $$m'(r_{n-1} - r) \leq s \leq \frac{p}{p-1}.$$  Then there exists $G' \in \G_n$, a polynomial
$H \in 1 + T^{-1}R[T^{-1}]$, and $F \in \KK$ such that
$\frac{G'}{G}H^p = 1 - p^sF$ with $v_r(F) = 0$ and $[F]_r = f$.
\end{cor}

\proof Let $F'$ be a polynomial in $T^{-1}$ of the same degree as $f$ such that $v_r(F') = 0$ and $[F']_r = f$.    
Now, $v_{r_{n-1}}(p^sF') = s - \deg(f)(r_{n-1} - r)$, which is positive by our assumptions.  
Choose $\sigma$ such that $s - (r_{n-1} - r) < \sigma < \frac{p}{p-1}$.  Then Lemma \ref{improvelem}(ii)
yields $G'$ and $H$ such that  $$\frac{G'}{G}H^p \equiv 1 - p^sF' \pmod{p^{\sigma}, \tilde{T}^{-m'}}.$$ 
Lemma \ref{improvelem}(ii) also allows us to assume that $$v_{r_{n-1}}(\frac{G'}{G}H^p - 1) \geq s - \deg(f)(r_{n-1} - r).$$

If $F$ is such that $\frac{G'}{G} H^p = 1 - p^sF$, then it suffices to show that $v_r(F - F') > 0$.  
Write $p^sF - p^sF'$ as a power series $\sum_{j=0}^{\infty} \alpha_j \tilde{T}^{-j}$.  For $j < m'$, we have $v(\alpha_j) \geq \sigma$.
For $j \geq m'$, we have $v(\alpha_j) \geq s - \deg(f)(r_{n-1}-r)$.  In both cases, we have $v(\alpha_j) + j(r_{n-1} - r) > s$, which
shows that $v_r(p^sF - p^sF') > s$, and thus $v_r(F - F') > 0$.
\Endproof

\subsection{Specialization to the context of Theorem \ref{Athm1} and Lemma \ref{Alem3}}\label{Ssimple}
Maintain the notation of \S\ref{Sarb}.  In this section, we show that Assumptions \ref{ass1} and \ref{ass2} are satisfied in the situations
of \S\ref{Sproof}.  In particular, we want to show that there is a $g$ satisfying Assumption \ref{ass1} in the form
\begin{equation} \label{eq2.0}
      g=t^{-N_2}\prod_{i=1}^{N_1}\,(1-\bar{x}_it^{-1}).
\end{equation}
Such a solution (where $a_0 = -N_2$ and $a_i = 1$ for $i > 1$) is called a solution with \emph{simple type}.
In the situation where we only look for solutions with simple type, \eqref{eq1.7} can be rewritten as
\begin{equation} \label{eq2.1}
  \omega := dg/g -m\sum_{j=0}^{\nu}\,t^{-mp^j-1}dt = \frac{cdt}{t^{N+m_{n-1}+1}g},
\end{equation}
where $c \ne 0$.
Recall that $N=m_n-m_{n-1} \geq m_{n-1}(p-1)$ and that $m_n \not\equiv 0\pmod{p}$
if the inequality is strict.  

\begin{lem}\label{Lsimple}
Suppose $g$ is of simple type as in \eqref{eq2.0}.  Then $g$ satisfies \eqref{eq2.1} (and thus \eqref{eq1.7}) if and only if it satisfies Assumption
\ref{ass1}.
\end{lem}

\proof
The ``only if" direction is immediate from the right hand side of \eqref{eq2.1}.  To prove the ``if" direction, suppose that $g$ is of simple type and 
satisfies Assumption \ref{ass1}.  Then $\omega$ has a zero of order $m_{n-1} + N_1 - 1$ at $\infty$, simple poles at all the $N_1$ zeroes of $g$,
a pole of order $m_{n-1} + 1$ at $0$, and no other poles.  Thus, $\omega$ can have no other zeroes.  So the differential forms
on each side of \eqref{eq2.1} have the same divisor, and we can choose $c$ to get our desired equality.  
\Endproof

\begin{thm} \label{thm2}
  Assume that $m = 1$ (i.e.\ $(m_1, \ldots, m_{n-1}) = (1, p, \ldots, p^{n-2} = p^{\nu})$). 
  Then \eqref{eq2.1} has a unique solution $g\in k(t)$. This solution is
  of the form \eqref{eq2.0}, with pairwise distinct $\bar{x}_i\in k$, and
  where
  \[
     p^{\nu+1}-p^{\nu}-p<N_1\leq p^{\nu+1}-p^{\nu}, \quad N_1\equiv N\pmod{p}
  \]

\end{thm}
    
\proof We rewrite \eqref{eq2.1} as the nonhomogenous linear differential
equation
\begin{equation} \label{eq2.2}
  dg -g\cdot\big(\sum_{j=0}^{\nu}t^{-p^j-1}dt\big) = ct^{-N-p^{\nu}-1}dt
\end{equation}
and we first look for solutions in $k((t))$ of the form
$g=\sum_{i=-N}^\infty \alpha_it^i$. We obtain a system of linear equations in the
coefficients $\alpha_i$.  In degree $-N - p^{\nu} - 1$, we obtain $\alpha_{-N}=-c$, and for $i \geq 1$
we get a linear expression for $\alpha_{-N+i}$ in terms of
$\alpha_{-N},\ldots,\alpha_{-N+i-1}$. In other words,
\eqref{eq2.2} has a unique solution of the form $g=-ct^{-N}+\ldots\in
k((t))$.  We also observe that the linear equations become periodic in $i$, as
soon as $i\geq p^{\nu}$, which means that the coefficients of $g$ (which only take values in $\FF_p$) are also eventually periodic. 
This means that $g$ is actually a \emph{rational} function in $t$, i.e.\ $g\in k(t)$ (as if $P$ is the period of the coefficients of $g$, then
one can write an equation relating $g$ and $t^Pg$).  Now we can use
\eqref{eq2.1} to see that
\begin{equation} \label{eq2.2.1}
     \ord_0(g)=-N,\quad \ord_{\bar{x}}(g)\in\{0,1\},\;\bar{x}\neq
     0,\infty. 
\end{equation}
It follows that $g\in k[t,t^{-1}]$ is a Laurent polynomial with leading term
$-ct^{-N}$ with only simple zeroes outside of $t=0,\infty$. Set
\begin{equation}\label{En2}
     N_2:=\ord_\infty(g)
\end{equation}
and $N_1:=N-N_2$. Then $g$ is of the form \eqref{eq2.0}. It remains to see
that $N_1$ is as claimed in the theorem.

Let $\omega$ denote the differential form in \eqref{eq2.1}. By \eqref{En2}
we have
\begin{equation} \label{eq2.5}
     \ord_\infty(\omega) = p^{\nu}-1+N_1\geq p^{\nu}-1\geq 0,
\end{equation}
which implies that ${\rm res}_\infty(\omega)= {\rm res}_\infty(dg/g)= 0$. We conclude that $N_2\equiv
0\pmod{p}$ and hence $N_1\equiv N\pmod{p}$. It remains to prove the inequality
for $N_1$ stated in the theorem. 

First assume $\nu = 0$.  Then $\ord_{\infty}(\omega) = N_1$. 
Also, $\C(\omega) = dg/g = \omega + t^{-2}dt$.  Thus, either $N_1 = 0$ or
$\ord_{\infty}\C(\omega)=0$.  But it is easy to see that if $\ord_{\infty}\C(\omega) = 0$ and 
$\ord_{\infty}(\omega) \geq 0$, then $\ord_{\infty}\omega < p$.  Thus $0 \leq N_1 < p$, which 
is exactly the condition of the theorem.

Now assume $\nu\geq 1$. Set 
\[
    h:=t^{-N_2}g^{-1}=\prod_{i=1}^{N_1}(1-\bar{x}_it^{-1})^{-1}.
\]
Since $N_2\equiv 0\pmod{p}$, $h$ satisfies the differential equation
\begin{equation} \label{eq2.9}
     dh + h\cdot\big(\sum_{j=0}^{\nu}t^{-p^j-1}dt\big) 
        = ct^{-N_1-p^{\nu}-1}h^2\,dt
\end{equation}
which is derived from \eqref{eq2.2}.  If we write
$h=1+b_1t^{-1}+b_2t^{-2}+\ldots$ and plug this into \eqref{eq2.9}, we see that
the coefficients $b_1,\ldots,b_{p-1}$ satisfy the simple linear equations
\[
     b_1=1,\;2b_2=b_1,\ldots,(p-1)b_{p-1}=b_{p-2}.
\]
We conclude that $b_i=1/i!\neq 0$ for $i=0,1,\ldots,p-1$. Write $N_1=pM-a$ with
$0\leq a<p$ and $M\in\ZZ$. Then
\[
  \omega = ct^{-N_1-p^{\nu}-1}h\,dt
    = c(t^{-M-p^{\nu-1}})^p(t^a+b_1t^{a-1}+\ldots)\,dt/t,
\]
from which we see that
\[
  \C(\omega) = c^{1/p} \cdot t^{-M-p^{\nu-1}}(b_a^{1/p}+b_{p+a}^{1/p}t^{-1}+\ldots)\,dt/t.
\]
Since $b_a\neq 0$, it follows that
\begin{equation} \label{eq2.12}
  \ord_\infty\C(\omega) = M+p^{\nu-1}-1.
\end{equation}
In particular, this and \eqref{eq2.5} show that
$\ord_\infty\C(\omega)<\ord_\infty\omega$. Therefore, the equality
\[
  \C(\omega) = \omega+t^{-p^{\nu}-1}dt.
\]
shows that
\begin{equation} \label{eq2.14}
  \ord_\infty\C(\omega) = p^{\nu}-1.
\end{equation}
Combining \eqref{eq2.12} and \eqref{eq2.14} we conclude that
$M=p^{\nu}-p^{\nu-1}$ and hence 
\[
      N_1=pM-a=p^{\nu+1}-p^{\nu}-a.
\]
This completes the proof of the theorem.
\Endproof

\begin{cor}\label{1tom}
If $m | N$, then (\ref{eq2.1}) has a unique solution $g(t) \in k(t)$.  In this case, the solution is given by
$g(t) = h(t^m)$, where $h(t)$ is the solution to (\ref{eq2.2}) with $N$ replaced by $N' := N/m$ and $c$ replaced by $c' := c/m$.  
In particular, $g(t) \in k[t^{-1}]$ 
is of simple type, and of degree $N$ in $t^{-1}$.  Furthermore, $m | N_1$ and we have
$$m_{n-1}(p-1)-mp<N_1\leq m_{n-1}(p-1), \quad N_1\equiv N\pmod{p}.$$
\end{cor}

\proof Equation (\ref{eq2.1}) can be rewritten as
\begin{equation} \label{eqwithm1}
dg - g \cdot m \cdot\big(\sum_{j=0}^{\nu}t^{-mp^j-1}dt\big) = ct^{-N-m_{n-1}-1}dt,
\end{equation}
Since $g(t) = h(t^m)$, we have $dg = mt^{m-1}dh(t^m)$.  By definition, $h$ satisfies
\begin{equation}\label{heq}
dh - h \cdot\big(\sum_{j=0}^{\nu}t^{-p^j-1}dt\big) = c't^{-N'-p^{\nu}-1}dt.
\end{equation}
Multiplying \eqref{heq} by $t$, substituting $t^m$ for $t$ on both sides, and then multiplying by $\frac{m}{t}$ yields (\ref{eq2.1}).  
This shows that $g(t)$ is as in the corollary.  Then $g$ is unique because its Laurent series coefficients can be calculated by recursion.  
The properties of $N_1$ follow easily. 
\Endproof

\begin{rem}\label{recursion}
Suppose $g \in k[t^{-1}]$ is a solution to (\ref{eq2.1}), thus (\ref{eqwithm1}).  
Let $\cb_i$ be the coefficient of $t^{-i}$ in $g$ (with $\cb_i = 0$ if $i < 0$).  Then 
(\ref{eqwithm1}) in degree $-i-1$ gives the equation
\begin{equation}\label{eqrec}
i\cb_i + m(\cb_{i-m} + \cb_{i - pm} + \cdots + \cb_{i - p^{\nu}m}) = \begin{cases} -c & i = N + m_{n-1} \\ 0 & \text{otherwise.} \end{cases}
\end{equation}
\end{rem}

Corollary \ref{1tom} proves Theorem \ref{Athm1}. 
Since $m | m_{n-1}$ by definition and $N = m_n - m_{n-1}$, we have that $m_n = pm_{n-1}$ implies $m | N$,
and thus that the solution $g$ to \eqref{eq2.1} guaranteed by Corollary \ref{1tom} satisfies Assumption \ref{ass1}.

\begin{lem}\label{Lassumptions}
If $g$ is a solution to \eqref{eq2.1} with simple type, and we further assume that $m_n = pm_{n-1}$,  
then Assumption \ref{ass2} holds as well.
\end{lem}

\proof  Corollary \ref{1tom} shows that $N_1 = m_{n-1}(p-1)$.  One then sees that the matrix $A$ in \eqref{eq3.6} is square.  
We show that $\ker(A) = 0$.  Suppose $\vec{v} \in \ker(A)$ is non-zero, and let $\vec{x}$ be the vector $(\xb_1, \ldots, \xb_n)$. 
Then, if $\epsilon^2 = 0$, the vector $\vec{x} + \epsilon\vec{v}$ must also satisfy Assumption \ref{ass1}, and thus, by Lemma \ref{Lsimple}, 
$\vec{x} + \epsilon\vec{v}$ must satisfy 
\eqref{eqwithm1} (equivalent to \eqref{eq2.1}).  Now, we claim that replacing $\vec{x}$ with $\vec{x} + \epsilon\vec{v}$ 
replaces $g$ with $g + \epsilon h(\vec{v})$, where $h(\vec{v})$ is a nonzero polynomial in $t^{-1}$.  Given the claim, 
Equation \eqref{eqwithm1} then implies that $dh - h \cdot m \cdot (\sum_{j=0}^{\nu} t^{-mp^j-1}dt) = 0$, or that
$$dh/h = m (\sum_{j=0}^{\nu} t^{-mp^j-1}dt).$$  But the right-hand side is not logarithmic, as it has a multiple pole at $t=0$.
This is a contradiction.

To prove the claim, we must show that $h(\vec{v})$ is nonzero.  If $\vec{e}_j$ is the $j$th standard basis vector, then
$h(\vec{e}_j) = t^{-N_2-1}\left(\prod_{i=1}^{N_1}\,(1-\xb_it^{-1})\right) / (1-\xb_jt^{-1})$.  In particular, $\xb_i$ is a zero of 
$h(\vec{e}_j)$ for all $i \ne j$, but not when $i=j$.  But then no $h(\vec{e}_j)$ can be a linear combination of the $h(\vec{e}_i)$ with $i \ne j$, because
that would imply that $\xb_j$ is a zero of $h(\vec{e}_j)$.  Thus $h(\vec{v}) \ne 0$ for $\vec{v} \ne 0$.
\Endproof

\begin{cor}\label{proof}
Lemma \ref{Alem3} holds.
\end{cor}

\proof If $m_n = pm_{n-1}$, then Corollary \ref{1tom} and Lemma \ref{Lassumptions} guarantee that there is a unique $g$ satisfying \eqref{eq2.1}
(thus Assumption \ref{ass1}) and Assumption \ref{ass2}, and that $N_2 = 0$.  
By Lemma \ref{Lsimple}, the family $\G_n$ from \S\ref{Sarb} (in particular, Corollary \ref{Capprox}) 
is the same as the family defined in Definition
\ref{Adef1}.  Then Lemma \ref{Alem3} follows from Corollary \ref{Capprox}, noting that, since $N = N_1 = m_{n-1}(p-1)$, we have
$m' = m_n = pm_{n-1}$ in Corollary \ref{Capprox}.
\Endproof

\begin{rem}\label{samebounds}
Suppose $m = 1$, so that $m_{n-1} = p^{\nu-2}$.  
Then the condition $p^{n-1} - p^{n-2} - p < N_1 \leq p^{n-1} - p^{n-2}$ from Theorem \ref{thm2} is exactly the same as the
necessary condition for the matrix $A$ in \eqref{eq3.6} to be square (Remark \ref{squarerem}), and thus it is necessary in order for
Assumption \ref{ass2} to be satisfied.   So if $g$ satisfies Assumption \ref{ass1}, even if $g$ does not have simple type, we still must have
the same bounds on $N_1$ in order to proceed with any proof that uses Assumption \ref{ass2}, which is essential for
our proof of Lemma \ref{Alem3}.  The importance of bounds on $N_1$ was shown in Example \ref{exa2} and Remark \ref{Rarb}. 
\end{rem}

\subsection{An invertible matrix}\label{Smatrix}
We give a result (Lemma \ref{lem3}) that will be used in \S\ref{Sbeyond}.
Assume that $m_n = pm_{n-1}$, and that $g$ is a simple type solution to \eqref{eq2.1}.  
By Corollary \ref{1tom}, we have $N_1 = N = m_{n-1}(p-1)$, and $N_2 = 0$.
Maintain the notation $\G_n$ from \S\ref{Sarb}.  Recall that $\tilde{T} = p^{-r_{n-1}}T$.  Recall also that, if $\Gamma \in \G_n$, then 
$\Gamma \in 1 + \tilde{T}^{-1}R[\tilde{T}^{-1}]$ has degree $m_{n-1}(p-1)$
and simple zeroes.

Let $V \cong R^{m_{n-1}}$ be the free $R$-module of polynomials of the form  
$$\sum_{j=1}^{m_{n-1}}b_j \tilde{T}^{-pj} \ \ (b_j \in R),$$ and $W \cong R^{m_{n-1}}$ be the free $R$-module of polynomials of the form
$$\sum_{j = m_{n-1}(p-1) + 1}^{pm_{n-1}} a_j \tilde{T}^{-j} \ \ (a_j \in R).$$ 
For $\Gamma \in \G_n$, consider the linear map $L:V \to W$ defined by
$$L_{\Gamma}(J) = (\Gamma J)^{tr},$$
where $(\cdot)^{tr}$ truncates a polynomial in $\tilde{T}^{-1}$ to preserve only the terms of degree $m_{n-1}(p-1) + 1$ through
$pm_{n-1}$, inclusive. 
Let $A_{\Gamma} \in M_{m_{n-1}}(R)$ be the matrix of this linear map, relative to the input basis 
$\tilde{T}^{-p}, \ldots, \tilde{T}^{-pm_{n-1}}$, and the output basis
$\tilde{T}^{-(m_{n-1}(p-1) + 1)}, \ldots, \tilde{T}^{-pm_{n-1}}$.  If $c_i$ is the coefficient of $\tilde{T}^{-i}$ in $\Gamma$ (with $c_i = 0$ for $i < 0$), 
then the matrix $A_{\Gamma}$ is given by
\begin{equation}\label{matrixeq}
A_{\Gamma} = \left( \begin{array}{cccc} 
c_{m_{n-1}(p-1) - p + 1} & c_{m_{n-1}(p-1) - 2p + 1} & \cdots & c_{-m_{n-1} + 1} \\
c_{m_{n-1}(p-1) - p + 2} & c_{m_{n-1}(p-1) - 2p + 2 }& \cdots & c_{-m_{n-1} + 2} \\
\vdots & \vdots & \ddots & \vdots \\
c_{pm_{n-1} - p} & c_{pm_{n-1} - 2p} & \cdots & c_0
\end{array} \right)
\end{equation}

\begin{lem}\label{lem2}
For any $\Gamma \in \G_n$, the matrix $A_{\Gamma}$ above lies in $GL_{m_{n-1}}(R)$.
\end{lem}

\proof We note that $A_{\Gamma}$ lying in  $GL_{m_{n-1}}(R)$ is equivalent to the reduction $\Ab_{\Gamma}$ of $A_{\Gamma}$ lying in 
$GL_{m_{n-1}}(k)$, which is
equivalent to the linear transformation given by $\Ab_{\Gamma}$ being surjective.  Since $\Ab_{\Gamma}$ does not depend on the choice of 
${\Gamma} \in \G_n$, it suffices to show that, for any $w \in W$, there exists \emph{some} $\Gamma \in \G_n$ and some
$v \in V$ such that $L_{\Gamma}(v) = w$. 
Consequently, it suffices exhibit an $R$-basis $w_1, \ldots, w_{m_{n-1}}$ of $W$ such that each $w_i$ is in the image of some
$L_{\Gamma}$.  

Take any $G \in \G_n$, and let $$w_i' = (\tilde{T}^{-(m_{n-1}(p-1) + i)}G)^{tr}, \quad  1 \leq i \leq m_{n-1}.$$  
One easily sees that the $w_i'$ form an $R$-basis of $W$.  
If $i = m_{n-1}$, take $w_{m_{n-1}} = w_{m_{n-1}}' = \tilde{T}^{-pm_{n-1}}$.  Then, $w_i = L_{\Gamma}(\tilde{T}^{-pm_{n-1}}),$
regardless of $\Gamma$.

Now, suppose $1 \leq i < m_{n-1}$.  For any small $\epsilon \in \QQ_{> 0}$, 
let $r \in \QQ$ be such that $$(m_{n-1} - i)(r_{n-1} - r) = \epsilon$$ (we choose
$\epsilon$ small enough so that $r > 0$).  If $s = pm_{n-1}(r_{n-1} - r)$, then the monomial 
$$F = p^{(m_{n-1}(p-1) + i)(r - r_{n-1})}\tilde{T}^{-(m_{n-1}(p-1) + i)}$$ satisfies 
$[F]_r = t^{-(m_{n-1}(p-1) + i)}$
and $(p^sFG)^{tr} = p^{\epsilon}w_i'$.
By Lemma \ref{improvelem}(i) with $J = 1 + p^sF$, there is a unique solution $(G', I)$ to $$G'I - G \equiv p^sFG \pmod{\tilde{T}^{-pm_{n-1}}},$$ 
where $G' \in \G_n$ and $I$ is of the form $1 + \sum_{j=1}^{m_{n-1} - 1}b_j \tilde{T}^{-pj}$, with all $b_j \in \m$.  
Then $L_{G'}(I - 1) = p^{\epsilon}w_i' + cw_{m_{n-1}}$, with $c \in R$ and $I - 1 \in V$.  By continuity and the 
linearity of $L_{G'}$, it follows that there exists $c' \in R$ and $v \in V$ such that $L_{G'}(v) = w_i' + c'w_{m_{n-1}}$.
If $w_i = w_i' + c'w_{m_{n-1}}$, then $w_i$ is in the image of $L_{G'}$.  Since $w_1, \ldots, w_{m_{n-1}}$ form a basis of $W$, we are done.
\Endproof

\begin{lem}\label{lem3}
Let $\alpha = (m_n - p[\frac{m_n}{p}])$.
Let $C \in M_{m_{n-1}}(R)$ be the matrix with 
$C_{ij} = c_{m_{n-1}(p-1)  - pj +i + \alpha}$.  Then $C \in GL_{m_{n-1}}(R)$.
\end{lem}

\proof It suffices to show that the reduction $\ol{C}$ of $C$ to $M_{m_{n-1}}(k)$ is invertible. 
For all $j$, let $\ol{c}_j$ be the reduction of $c_j$ to $k$.  
These are the coefficients of $t^{-j}$ in $g$ (see Remark \ref{recursion}).

By Lemma \ref{lem2}, the reduction
$\ol{A}$ of $A_{\Gamma}$ in (\ref{matrixeq}) to $M_{m_{n-1}}(k)$ is invertible. 
Notice that $\ol{A}$ satisfies $\ol{A}_{ij} = \ol{c}_{m_{n-1}(p-1) - pj + i}$. 
Thus, up to reordering rows, the matrix $\ol{C}$ is derived from the matrix $\ol{A}$ by replacing each of the first $\alpha$ rows 
$R_1, \ldots, R_{\alpha}$ of $\ol{A}$ with rows $R_1', \ldots, R_{\alpha}'$, such that each entry $\ol{c}_\ell$ in $R_i$ is replaced 
by $\ol{c}_{\ell + m_{n-1}}$ in $R_i'$.

Now, setting $i = l + m_{n-1}$ in (\ref{eqrec}), we obtain that the $c_l$ satisfy the recursion
\begin{equation}\label{Erec}
\ol{c}_{\ell+m_{n-1}} = -\frac{m}{\ell+m_{n-1}}(\ol{c}_{\ell + m_{n-1} - m} + \ol{c}_{\ell + m_{n-1} - pm} + \cdots + \ol{c}_{\ell + m_{n-1} - p^{\nu}m}),
\end{equation} so long as $p \nmid \ell + m_{n-1}$. 
Since $0 \leq \alpha < p$, we see that in the first $\alpha$ rows of $\ol{A}$, no $\ol{c}_l$ appears with $p|\ell+m_{n-1}$, 
so (\ref{Erec}) holds.  All entries of $R_1$ are of the form $\ol{c}_\ell$, with $\ell \equiv 1 - m_{n-1} \pmod{p}$, so (\ref{Erec}) shows that $R_1'$ is a 
linear combination of the rows of $\ol{A}$, where the coefficient of $R_1$ is the unit
$-\frac{m}{1}$.  So replacing $R_1$ with $R_1'$ gives a matrix $\ol{A}_1$ such that
$\det(\ol{A}_1) = -m\det(\ol{A})$.  In particular, $\ol{A}_1$ is invertible.  For the same reasons, replacing the row $R_2$ of
$\ol{A}_1$ with $R_2'$ gives a matrix $\ol{A}_2$ which is again invertible.  Repeating this process a total of $\alpha$ times yields an invertible
matrix $\ol{A}_{\alpha}$.  Since $\ol{C}$ is obtained from $\ol{A}_{\alpha}$ by reordering rows, it is also invertible.

\Endproof

\subsection{Proof of Propositions \ref{Pnonminimal} and \ref{Pwithindisk}}\label{Sbeyond}
Maintain the notation and assumptions of \S\ref{Smatrix}.  Let $G_{\min} \in \G_n$.  Note that $G_{\min} \in R[\tilde{T}^{-1}]$ has degree 
$m_n - m_{n-1} = m_{n-1}(p-1)$.  Let $m_n' \geq m_n = pm_{n-1}$, with $m_n'$ prime to $p$ unless $m_n' = m_n$. 
Let $F \in T^{-1}R[T^{-1}]$ have degree $\leq m_n'$.

In order to prove Proposition \ref{Pnonminimal}, we begin by trying to find polynomials $I$ and $G_n'$ in $R[T]$ such that 
\begin{equation}\label{Emodify}
G_{\min}I - G_n' \equiv -p^{\frac{p}{p-1}}F \pmod{p^{\frac{p}{p-1} + \epsilon}R[T^{-1}]}
\end{equation} 
for some $\epsilon > 0$, with $I$ having only monomials of degree divisible by $p$ and $G_n'$ having degree at most $m_n' - m_{n-1}$.
Consider instead the congruence
\begin{equation}\label{Emodify2}
G_{\min}I - G_n' \equiv -p^{\frac{p}{p-1}}F \pmod{\tilde{T}^{-m_n'-1}K[\tilde{T}^{-1}]}.
\end{equation}
Since $F$ has degree at most $m_n'$ in $\tilde{T}^{-1}$, we look for $I$ of the form 
$$I = \sum_{i=0}^{[\frac{m_n'}{p}]} b_i\tilde{T}^{-pi}.$$ 
Being able to choose $G_n'$ affords us some freedom; we   
may choose the coefficients $b_i$ for $i \leq \frac{m_n'}{p} - m_{n-1}$ at will.  This is because the terms in $G_{\min}I$ 
to which these $b_i$ contribute have degree at most $m_n' - m_{n-1}$.
As in \S\ref{Smatrix}, write $G_{\min} = \sum_{i=0}^{m_{n-1}(p-1)} c_i \tilde{T}^{-i}$, with $c_i = 0$ for $i$ not between $0$ and $m_{n-1}(p-1)$.

\begin{prop}\label{Pmodifiedsolution}
After a possible finite extension of $K$, there exists a solution $(I, G_n')$ 
to the equation $G_{\min}I - G_n' \equiv -p^{\frac{p}{p-1}}F \pmod{\tilde{T}^{-m_n' - 1}K[\tilde{T}^{-1}]}$ 
as described above with coefficients in $K$.  Furthermore,
a solution exists with $b_0 = 1$, $b_i = 0$ for $1 \leq i \leq \frac{m_n'}{p} - m_{n-1}$, and $v(b_i) \geq \frac{p}{p-1} - \frac{m_n'}{m_{n-1}(p-1)}$ for all other $i$.  
\end{prop}

\proof
We take $b_0 = 1$ and $b_i = 0$ for $1 \leq i \leq \frac{m_n'}{p} - m_{n-1}$.   
We rewrite (\ref{Emodify}) as $$G_{\min}(I-1) - G_n' \equiv p^{\frac{p}{p-1}}F - G_{\min} \pmod{\tilde{T}^{-m_n'-1}K[\tilde{T}^{-1}]}.$$  
This is a system of $m_n' + 1$ equations 
(corresponding to the terms of degree $0$ through $m_n'$) in $m_n' + 1$ variables 
(corresponding to the coefficients $c_0$ through $c_{m_n' - m_{n-1}}$ of $G_n'$ and 
$b_{[\frac{m_n'}{p}] - m_{n-1} + 1}$ through $b_{[\frac{m_n'}{p}]}$ of $I$).  
Let $E$ be the matrix representing this system.  We will show that $E \in GL_{m_n'+1}(R)$.  
Let the columns of $E$ correspond to the variables just listed, and let the rows of $E$ correspond to degrees $0$ through $m_n'$, in order.  
Then $E$ can be written as a block matrix,

\begin{equation}\label{blockeq}
E = \blockmatrix{-I_{m_n' - m_{n-1} + 1}}{B}{0}{C},
\end{equation}
where $I_{m_n' - m_{n-1} + 1}$ is the identity matrix of size $m_n' - m_{n-1} + 1$. 
Since the entries of $B$ and $C$ lie in $R$, being coefficients of $G_{\min}$, we see that showing that $E \in GL_{m_n'+1}(R)$ is 
equivalent to showing that $C \in GL_{m_{n-1}}(R)$.  Now, the $ij$th 
entry of $C$ is $c_\ell$, with $\ell = m_n' - m_{n-1} + i - p([\frac{m_n'}{p}]-m_{n-1} + j) = m_{n-1}(p-1) -pj + i + (m_n' - p[\frac{m_n'}{p}])$.  By
Lemma \ref{lem3}, the matrix $C$ lies in $GL_{m_{n-1}}(R)$, thus so does $E$.

Lastly, we show that we have $v(b_i) \geq \frac{p}{p-1} - \frac{m_n'}{m_{n-1}(p-1)}$ for all $i > \frac{m_n'}{p} - m_{n-1}$.  
Since $F \in R[T^{-1}]$, the coefficient of $\tilde{T}^{-i}$ in $-p^{\frac{p}{p-1}}F$ has valuation at least $\frac{p}{p-1} - \frac{i}{m_{n-1}(p-1)}$.  
Since $F$ has degree at most
$m_n'$ in $\tilde{T}^{-1}$, all coefficients of $-p^{\frac{p}{p-1}}F$ (in terms of $\tilde{T}$) have valuation at least 
$\frac{p}{p-1} - \frac{m_n'}{m_{n-1}(p-1)}$.  The $b_i$ in question are given as 
entries of the column vector $E^{-1}\vec{v}$, where $\vec{v}$ is the column vector whose entries are the coefficients of $1, \tilde{T}^{-1}, \ldots, 
\tilde{T}^{-m_n'}$ in $-p^{\frac{p}{p-1}}F$. Since $E^{-1}$ has entries in $R$, we are done.
\Endproof

We now prove Proposition \ref{Pnonminimal}.
\begin{cor}\label{Cmodifiedsolution}
Let $$I = \sum_{i=0}^{[\frac{m_n'}{p}]} b_i \tilde{T}^{-pi}$$ be the solution to (\ref{Emodify2}) found in Proposition \ref{Pmodifiedsolution}.
If $$H = \sum_{i=0}^{[\frac{m_n'}{p} ]} b_i^{1/p}\, \tilde{T}^{-i},$$ for any choice of $p$th roots, then there exists (after
a possible finite extension of $K$)
$G_n \in R[T^{-1}]$ of degree at most $m_n' - m_{n-1}$ such that $(H, G_n)$ is a solution to 
$$G_{\min}H^p - G_n \equiv -p^{\frac{p}{p-1}}F \pmod{p^{\frac{p}{p-1} + \epsilon}R[T^{-1}]}$$ 
for some $\epsilon \in \QQ_{>0}$.  Furthermore, $H$ and $G_n$ lie in $1 + T^{-1}\m[T^{-1}]$.
\end{cor}

\proof
Working in terms of $\tilde{T}$, we have that $H^p - I = \sum_{i = 1}^{p[\frac{m_n'}{p}]} pa_i \tilde{T}^{-i}$, where for all $i$,  
$v(a_i) \geq \min_j(v(b_j)) \geq \frac{p}{p-1} - \frac{m_n'}{m_{n-1}(p-1)}$.  Since $G_{\min} \in R[\tilde{T}^{-1}]$, it follows that
$$G_{\min}(H^p - I) \equiv \sum_{i=1}^{m_n'} e_i \tilde{T}^{-i} \pmod{\tilde{T}^{-m_n' - 1}R[\tilde{T}^{-1}]},$$ where 
$v(e_i) \geq 1 + \frac{p}{p-1} - \frac{m_n'}{m_{n-1}(p-1)}$ for all $i$.  Set $G_n = G_n' + \sum_{i=1}^{m_n' - m_{n-1}} e_i\tilde{T}^{-i}$, where $G_n'$ is 
as in Proposition \ref{Pmodifiedsolution}.  We must show that
$G_{\min}H^p - G_n \equiv -p^{\frac{p}{p-1}}F \pmod{p^{\frac{p}{p-1} + \epsilon}R[T^{-1}]}$ for small enough $\epsilon \in \QQ_{>0}$

If $0 \leq i \leq m_n' - m_{n-1}$, then the terms involving $\tilde{T}^{-i}$ (and thus $T^{-i}$)  
agree exactly for $G_{\min}H^p - G_n$ and $-p^{\frac{p}{p-1}}F$ by construction.  

If $m_n' - m_{n-1} +1 \leq i \leq m_n'$, then the terms involving $\tilde{T}^{-i}$ (or $T^{-i}$) agree exactly for $G_{\min}I$ and 
$-p^{\frac{p}{p-1}}F$ ($G_n$ and $G_n'$ have 
no terms of these degrees).  So the coefficient of $\tilde{T}^{-i}$ in $G_{\min}H^p - G_n + \frac{p}{p-1}F$ is $e_i$, which has valuation at least
$1 + \frac{p}{p-1} - \frac{m_n'}{m_{n-1}(p-1)}$.  Thus the coefficient of $T^{-i}$ in $G_{\min}H^{p} - G_n + p^{\frac{p}{p-1}}F$ has valuation at least
$$1 + \frac{p}{p-1} - \frac{m_n'}{m_{n-1}(p-1)} + \frac{i}{m_{n-1}(p-1)} > \frac{p}{p-1} + \frac{p-2}{p-1} \geq \frac{p}{p-1},$$ so the desired congruence
holds for these terms.  

For $i > m_n'$, the coefficient of $\tilde{T}^{-i}$ in $G_{\min}H^p - G_n + p^{\frac{p}{p-1}}F$ 
(which is the same as in $G_{\min}H^p$, as $F$ is
of degree $\leq m_n'$) has valuation at least $\frac{p}{p-1} - \frac{m_n'}{m_{n-1}(p-1)}$.  
So the corresponding coefficient of $T^{-i}$ has valuation at least
$$\frac{p}{p-1} - \frac{m_n'}{m_{n-1}(p-1)} + \frac{i}{m_{n-1}(p-1)} > \frac{p}{p-1},$$ proving the desired congruence.

To prove the last assertion, note that any non-constant coefficient of
$H$ is of the form $b_i^{1/p}\tilde{T}^i$, with $v(b_i) \geq \frac{p}{p-1} - \frac{m_n'}{m_{n-1}(p-1)}$ and $i > \frac{m_n'}{p} - m_{n-1}$.  It follows
easily that $v(b_i^{1/p}) + ir_{n-1} > 0$.  Thus, when $H$ is written in terms of $T^{-1}$, all non-constant coefficients have positive valuation.
The same is then true for $G_n$.
\Endproof

We examine our solution above in greater detail to prove Proposition \ref{Pwithindisk}.
\begin{prop}\label{Pconditionaldisk}
Let $G_n$ and $H$ be as in Corollary \ref{Cmodifiedsolution}.  Suppose that there is no integer $a$ satisfying
\begin{equation}\label{aineq}
\frac{m_n'}{p} - m_{n-1} < a \leq \left(\frac{m_n'}{m_n' - m_{n-1}}\right)\left(\frac{m_n'}{p} - m_{n-1}\right).
\end{equation} 
Then $v_{r_n'}(H - 1) > 0$ and $v_{r_n'}(G_n-1) > 0$, where $r_n' = \frac{1}{m_n'(p-1)}$.  
Thus all zeroes of $G_n$ lie in the open disk $D(r_n')$.
\end{prop}

\proof
Since $v_{r_n'}(G_{\min} - 1) > 0$, we need only show that $v_{r_n'}(H - 1) > 0$.    
Let $T' = p^{-r_n'}T$.  Writing $H - 1$ as $\sum_{i} \gamma_i (T')^{-i}$, we see from
Proposition \ref{Pmodifiedsolution} and the definition of $H$ 
that $\gamma_i$ is nonzero only when $[\frac{m_n'}{p}] - m_{n-1} + 1 \leq i \leq 
[\frac{m_n'}{p}]$.  In particular, $i > \frac{m_n'}{p} - m_{n-1}$.  Furthermore, if $b_i$ is the coefficient of $\tilde{T}^{-pi}$ in $I$ (Proposition
\ref{Pmodifiedsolution}), then by the definition of $H$ (Corollary \ref{Cmodifiedsolution}), we have
\begin{eqnarray*}
v(\gamma_i) & =& \frac{v(b_i)}{p} + i\left(\frac{1}{m_{n-1}(p-1)} - \frac{1}{m_n'(p-1)}\right) \\
&\geq& \frac{1}{p-1}\left(1 - \frac{m_n'}{pm_{n-1}} + \frac{i}{m_{n-1}} - \frac{i}{m_n'}\right)  
\end{eqnarray*}  
If $i > \left(\frac{m_n'}{m_n' - m_{n-1}}\right)(\frac{m_n'}{p} - m_{n-1})$, then $v(\gamma_i) > 0$, and we are done. 
\Endproof

\begin{exa}\label{Econditionaldisk}
We give a counterexample to Proposition \ref{Pconditionaldisk} when there is an $a$ satisfying \eqref{aineq}.
Let $p=5$, $n=3$, $(m_1, m_2, m_3) = (1, 5, 34)$, and $F = T^{-34}$.  Take $G_{\min} \in \G_3$ as above.
Then the matrix $C$ from \eqref{blockeq} becomes
$$\left( \begin{array}{ccccc} 
c_{20} & c_{15} & c_{10} & c_5 & c_0 \\
c_{21} & c_{16} & c_{11} & c_6 & c_1 \\
c_{22} &  c_{17} & c_{12} & c_7 & c_2 \\
c_{23} & c_{18} & c_{13} & c_8 & c_3 \\
c_{24} & c_{19} & c_{14} & c_9 & c_4 \\ \end{array} \right)$$
Now, based on (\ref{eqrec}) and the fact that $c_0$ must be $1$, one can calculate that the reduction of this matrix to $M_5(k)$ is

$$\ol{C} = \left( \begin{array}{ccccc} 
4 & 0 & 0 & 0 & 1 \\
0 & 1 & 1 & 1 & 1 \\
0 & 0 & 4 & 2 & 3 \\
0 & 0 & 0 & 1 & 1 \\
0 & 0 & 0 & 0 & 4 \\
\end{array} \right)$$
and $\ol{C}^{-1}$ has a nonzero entry in the upper righthand corner (easily checked using the adjoint formula for the inverse).

Now, if $E$ is as in \eqref{blockeq}, then $E^{-1} = \blockmatrix{-I_{30}}{BC^{-1}}{0}{C^{-1}}$.  
We calculate the coefficient $b_2$ of $\tilde{T}^{-10}$ in $I$ (from Proposition \ref{Pmodifiedsolution}).  It is the entry of
$E^{-1}\vec{v}$ (see end of proof of Proposition \ref{Pmodifiedsolution}) corresponding to the top row 
$\vec{w}$ of $C^{-1}$ (i.e., the row $(0 | \vec{w})$ of $E^{-1}$).  Here $\vec{v}$ is a vector whose only nonzero entry is in the last position, and is the
coefficient of $\tilde{T}^{-34}$ in $5^{5/4}F$, which has valuation $\frac{5}{4} - 34(\frac{1}{5(4)}) = -\frac{9}{20}$.  
Since $C^{-1}$ has an entry of valuation zero in the upper righthand corner, we see that 
$v(b_2) = -\frac{9}{20}$.  If $T' = 5^{-r_3}T = 5^{r_{2} - r_{3}}\tilde{T}$, then if we write $I$ in terms of $T'$, 
the coefficient of $(T')^{-10}$ has valuation $-\frac{9}{20} + \frac{10}{20} - \frac{10}{34(4)} = \frac{1}{20} - \frac{10}{136} < 0$.  The coefficient
of $(T')^{-10}$ in $H^p$ has the same valuation, as does the coefficient of $(T')^{-10}$ in $G_{\min}H^p$, as does the coefficient of $(T')^{-10}$ in 
$G_3$.  So the zeroes of $G_3$ do not all lie in the disk $D(r_3)$.  In fact, 
using the Newton polygon, one can show that $10$ zeroes of $G_3$ have valuation $\frac{1}{200}$. 
\end{exa}

\bibliographystyle{hplain}
\bibliography{literatur}

\begin{thebibliography}{10}

\bibitem{ArzdorfWewers}
Kai Arzdorf and Stefan Wewers.
\newblock A local proof of the semistable reduction theorem.
\newblock preprint, 2011.

\bibitem{Bertin}
Jos{\'e} Bertin.
\newblock Obstructions locales au rel\`evement de rev\^etements galoisiens de
  courbes lisses.
\newblock {\em C. R. Acad. Sci. Paris S\'er. I Math.}, 326(1):55--58, 1998.

\bibitem{BM}
Jos{\'e} Bertin and Ariane M{\'e}zard.
\newblock D\'eformations formelles des rev\^etements sauvagement ramifi\'es de
  courbes alg\'ebriques.
\newblock {\em Invent. Math.}, 141(1):195--238, 2000.

\bibitem{BGR}
Siegfried Bosch, Ulrich G{\"u}ntzer, and Reinhold Remmert.
\newblock {\em Non-{A}rchimedean analysis}, volume 261 of {\em Grundlehren der
  Mathematischen Wissenschaften [Fundamental Principles of Mathematical
  Sciences]}.
\newblock Springer-Verlag, Berlin, 1984.
\newblock A systematic approach to rigid analytic geometry.

\bibitem{BLstable}
Siegfried Bosch and Werner L{\"u}tkebohmert.
\newblock Stable reduction and uniformization of abelian varieties. {II}.
\newblock {\em Invent. Math.}, 78(2):257--297, 1984.

\bibitem{BLrigid}
Siegfried Bosch and Werner L{\"u}tkebohmert.
\newblock Formal and rigid geometry. {I}. {R}igid spaces.
\newblock {\em Math. Ann.}, 295(2):291--317, 1993.

\bibitem{arizona}
Irene Bouw and Stefan Wewers.
\newblock Group actions on curves and the lifting problem.
\newblock Course notes for the Arizona Winter School, 2012, available at
  http://swc.math.arizona.edu/aws/2012.

\bibitem{Brewisthesis}
Louis Brewis.
\newblock Ramification theory of the $p$-adic open disc and the lifting
  problem.
\newblock Ph.D. Thesis, available at
  http://webdoc.sub.gwdg.de/ebook/dissts/Ulm/Brewis2009.pdf.

\bibitem{Brewis}
Louis~Hugo Brewis and Stefan Wewers.
\newblock Artin characters, {H}urwitz trees and the lifting problem.
\newblock {\em Math. Ann.}, 345(3):711--730, 2009.

\bibitem{CGH1}
Ted Chinburg, Robert Guralnick, and David Harbater.
\newblock Oort groups and lifting problems.
\newblock {\em Compos. Math.}, 144(4):849--866, 2008.

\bibitem{Epp73}
Helmut~P. Epp.
\newblock Eliminating wild ramification.
\newblock {\em Invent. Math.}, 19:235--249, 1973.

\bibitem{Garuti96}
Marco~A. Garuti.
\newblock Prolongement de rev\^etements galoisiens en g\'eom\'etrie rigide.
\newblock {\em Compositio Math.}, 104(3):305--331, 1996.

\bibitem{Garuti99}
Marco~A. Garuti.
\newblock Linear systems attached to cyclic inertia.
\newblock In {\em Arithmetic fundamental groups and noncommutative algebra
  ({B}erkeley, {CA}, 1999)}, volume~70 of {\em Proc. Sympos. Pure Math.}, pages
  377--386. Amer. Math. Soc., Providence, RI, 2002.

\bibitem{Gre04}
Barry Green.
\newblock Realizing deformations of curves using {L}ubin-{T}ate formal groups.
\newblock {\em Israel J. Math.}, 139:139--148, 2004.

\bibitem{GM98}
Barry Green and Michel Matignon.
\newblock Liftings of {G}alois covers of smooth curves.
\newblock {\em Compositio Math.}, 113(3):237--272, 1998.

\bibitem{GM99}
Barry Green and Michel Matignon.
\newblock Order {$p$} automorphisms of the open disc of a {$p$}-adic field.
\newblock {\em J. Amer. Math. Soc.}, 12(1):269--303, 1999.

\bibitem{Ha80}
David Harbater.
\newblock Moduli of {$p$}-covers of curves.
\newblock {\em Comm. Algebra}, 8(12):1095--1122, 1980.

\bibitem{KatoVC}
Kazuya Kato.
\newblock Vanishing cycles, ramification of valuations, and class field theory.
\newblock {\em Duke Math. J.}, 55(3):629--659, 1987.

\bibitem{Ka86}
Nicholas~M. Katz.
\newblock Local-to-global extensions of representations of fundamental groups.
\newblock {\em Ann. Inst. Fourier (Grenoble)}, 36(4):69--106, 1986.

\bibitem{Lang}
Serge Lang.
\newblock {\em Algebra}, volume 211 of {\em Graduate Texts in Mathematics}.
\newblock Springer-Verlag, New York, third edition, 2002.

\bibitem{LiuAG}
Qing Liu.
\newblock {\em Algebraic geometry and arithmetic curves}, volume~6 of {\em
  Oxford Graduate Texts in Mathematics}.
\newblock Oxford University Press, Oxford, 2002.
\newblock Translated from the French by Reinie Ern{\'e}, Oxford Science
  Publications.

\bibitem{Nagata}
Masayoshi Nagata.
\newblock {\em Local rings}.
\newblock Robert E. Krieger Publishing Co., Huntington, N.Y., 1975.
\newblock Corrected reprint.

\bibitem{Ob11}
Andrew Obus.
\newblock The (local) lifting problem for curves.
\newblock Preprint, 2011, arXiv:1105.1530.

\bibitem{Oort95}
Frans Oort.
\newblock Some questions in algebraic geometry.
\newblock Available at http://www.staff.science.uu.nl/$\sim$oort0109/.

\bibitem{Oort87}
Frans Oort.
\newblock Lifting algebraic curves, abelian varieties, and their endomorphisms
  to characteristic zero.
\newblock In {\em Algebraic geometry, {B}owdoin, 1985 ({B}runswick, {M}aine,
  1985)}, volume~46 of {\em Proc. Sympos. Pure Math.}, pages 165--195. Amer.
  Math. Soc., Providence, RI, 1987.

\bibitem{OSS}
Frans Oort, Tsutomi Sekiguchi, and Noriyuki Suwa.
\newblock On the deformation of {A}rtin-{S}chreier to {K}ummer.
\newblock {\em Ann. Sci. \'Ecole Norm. Sup. (4)}, 22(3):345--375, 1989.

\bibitem{Pop}
Florian Pop.
\newblock Lifting of curves: The {O}ort conjecture.
\newblock preprint, 2012, arXiv:1203.1867.

\bibitem{SS94}
Tsutomi Sekiguchi and Noriyuki Suwa.
\newblock On the unified {K}ummer-{A}rtin-{S}chreier-{W}itt theory.
\newblock Chuo Math. preprint series, 1994.

\bibitem{SS99}
Tsutomi Sekiguchi and Noriyuki Suwa.
\newblock On the unified {K}ummer-{A}rtin-{S}chreier-{W}itt theory.
\newblock Laboratoire de Math\'{e}matiques Pures de Bordeaux preprint series,
  1999.

\bibitem{SS01}
Tsutomu Sekiguchi and Noriyuki Suwa.
\newblock A note on extensions of algebraic and formal groups. {IV}.
  {K}ummer-{A}rtin-{S}chreier-{W}itt theory of degree {$p\sp 2$}.
\newblock {\em Tohoku Math. J. (2)}, 53(2):203--240, 2001.

\bibitem{SerreCL}
Jean-Pierre Serre.
\newblock {\em Corps locaux}.
\newblock Hermann, Paris, 1968.
\newblock Deuxi{\`e}me {\'e}dition, Publications de l'Universit{\'e} de
  Nancago, No. VIII.

\bibitem{Tossici10}
Dajano Tossici.
\newblock Models of {$\mu\sb {p\sp 2,K}$} over a discrete valuation ring.
\newblock {\em J. Algebra}, 323(7):1908--1957, 2010.
\newblock With an appendix by Xavier Caruso.

\bibitem{boundary}
Stefan Wewers.
\newblock Swan conductors on the boundary of {L}ubin-{T}ate spaces.
\newblock preprint, 2005, arXiv:math.NT/0511434.

\bibitem{cyclic}
Stefan Wewers.
\newblock Fiercely ramified cyclic extensions of $p$-adic fields with residue
  field of dimension one.
\newblock Preprint, 2010, arXiv:1104.3785v1.

\bibitem{Witt37}
Ernst Witt.
\newblock {Zyklische K\"orper und Algebren der Charakteristik $p$ vom Grad
  $p^n$.}
\newblock {\em J.\ Reine Angew.\ Math.}, 176:126--140, 1937.

\bibitem{ZariskiSamuelII}
Oscar Zariski and Pierre Samuel.
\newblock {\em Commutative algebra. {V}ol. {II}}.
\newblock Springer-Verlag, New York, 1975.
\newblock Reprint of the 1960 edition, Graduate Texts in Mathematics, Vol. 29.

\end{thebibliography}

\end{document}